\documentstyle[11pt,bezier,figure]{article}

\begin{document}

\font\bbbld=msbm10 scaled\magstep1
\newcommand{\bfR}{\hbox{\bbbld R}}
\newcommand{\bfC}{\hbox{\bbbld C}}
\newcommand{\bfZ}{\hbox{\bbbld Z}}
\newcommand{\bfH}{\hbox{\bbbld H}}
\newcommand{\bfQ}{\hbox{\bbbld Q}}
\newcommand{\bfN}{\hbox{\bbbld N}}
\newcommand{\bfP}{\hbox{\bbbld P}}
\newcommand{\bfT}{\hbox{\bbbld T}}
\def\Sym{\mathop{\rm Sym}}
\newcommand{\halo}[1]{\Int(#1)}
\def\Int{\mathop{\rm Int}}
\def\Re{\mathop{\rm Re}}
\def\Im{\mathop{\rm Im}}
\newcommand{\union}{\cup}
\newcommand{\goesto}{\rightarrow}
\newcommand{\bdy}{\partial}
\newcommand{\n}{\noindent}
\newcommand{\p}{\hspace*{\parindent}}

\newtheorem{theorem}{Theorem}[section]
\newtheorem{assertion}{Assertion}[section]
\newtheorem{proposition}{Proposition}[section]
\newtheorem{lemma}{Lemma}[section]
\newtheorem{definition}{Definition}[section]
\newtheorem{claim}{Claim}[section]
\newtheorem{corollary}{Corollary}[section]
\newtheorem{observation}{Observation}[section]
\newtheorem{conjecture}{Conjecture}[section]
\newtheorem{question}{Question}[section]
\newtheorem{example}{Example}[section]

\newbox\qedbox
\setbox\qedbox=\hbox{$\Box$}
\newenvironment{proof}{\smallskip\noindent{\bf Proof.}\hskip \labelsep}%
                        {\hfill\penalty10000\copy\qedbox\par\medskip}
\newenvironment{remark}{\smallskip\noindent{\bf Remark.}\hskip \labelsep}%
                        {\hfill\penalty10000\copy\qedbox\par\medskip}
\newenvironment{remark1}{\smallskip\noindent{\bf Remark 1.}\hskip \labelsep}%
                        {\hfill\penalty10000\copy\qedbox\par\medskip}
\newenvironment{remark2}{\smallskip\noindent{\bf Remark 2.}\hskip \labelsep}%
                        {\hfill\penalty10000\copy\qedbox\par\medskip}
\newenvironment{proofspec}[1]%
                      {\smallskip\noindent{\bf Proof of Theorem 1.1.}
                        \hskip \labelsep}%
                        {\nobreak\hfill\hfill\nobreak\copy\qedbox\par\medskip}
\newenvironment{proofspec2}[1]%
                      {\smallskip\noindent{\bf Proof of Theorem 1.2.}
                        \hskip \labelsep}%
                        {\nobreak\hfill\hfill\nobreak\copy\qedbox\par\medskip}
\newenvironment{acknowledgements}{\smallskip\noindent{\bf Acknowledgements.}%
        \hskip\labelsep}{}

\setlength{\baselineskip}{1.0\baselineskip}

\title{On the Index of Constant Mean Curvature 1 Surfaces 
in Hyperbolic Space}
\author{Levi Lopes  de Lima, Wayne Rossman
\vspace{0.2in}\\
Department of Mathematics\\
Universidade Federal do Cear\'{a}\\
Fortaleza, Brasil}

\maketitle

\begin{abstract}
We show that the index of a constant mean curvature 1 surface in 
hyperbolic 3-space is completely determined by the compact Riemann surface 
and secondary Gauss map that represent it in Bryant's Weierstrass 
representation.  We give three applications of this observation.
Firstly, it allows us to explicitly compute the 
index of the catenoid cousins and some other examples.  
Secondly, it allows us to be able to apply a method similar to that of 
Choe (using Killing vector fields on 
minimal surfaces in Euclidean 3-space) to our case as well, 
resulting in lower bounds of index for other examples.  And thirdly, it 
allows us to give a more direct proof of the result by do 
Carmo and Silveira 
that if a constant mean curvature 1 surface in hyperbolic 3-space has 
finite total curvature, then it has finite index.  
Finally, we show that for 
any constant mean curvature 1 surface in hyperbolic 3-space that 
has been constructed via a correspondence to a minimal surface 
in Euclidean 3-space, we can take advantage of this correspondence 
to find a lower bound for its index.  
\end{abstract}

\section{Introduction}

In a seminal paper \cite{By}, R. Bryant has shown that the geometry of 
surfaces with constant mean curvature 1 in hyperbolic 3-space 
$\bfH^3(-1)$ has 
many similarities with the geometry of minimal surfaces in Euclidean 
space $\bfR^3$.  It was shown in particular that such surfaces 
admit a Weierstrass representation in terms of certain holomorphic 
data (see section 3 below for details).  A detailed analysis of this 
representation has allowed the construction of many complete examples 
(\cite{UY1}, \cite{RUY}).

It is well known that constant mean curvature surfaces 
in $\bfR^3$ and $\bfH^3$ can be 
characterized as critical points for the area functional, under 
compactly supported variations.  (Recall that in the constant mean 
curvature nonzero case, only volume preserving variations are allowed.)  
Regarding these variational problems, in the cases of complete minimal 
surfaces in $\bfR^3$ and complete constant mean curvature 1 surfaces in 
$\bfH^3$, it is known that the only stable 
objects are planes and horospheres \cite{CP}, \cite{Si}.  
This makes all the more 
interesting the study of surfaces of finite index, namely, surfaces 
for which the dimension of the space of area decreasing variations is 
finite.  A fundamental result regarding this point is due to 
Fischer-Colbrie \cite{FC}, who has shown that a minimal surface 
in $\bfR^3$ has finite index if and only if its total curvature is 
finite.  (In regard to this, see also \cite{G}.)  In fact, 
Fischer-Colbrie's analysis allows us to obtain 
explicit estimates for the index of concrete examples with finite 
total curvature.   
Recall that any such surface is conformally equivalent to a compact 
Riemann surface $\Sigma$ punctured at finitely many points 
corresponding to the ends of the original surface.  Moreover, the Gauss 
map of $M$ extends meromorphically across the punctures defining a 
meromorphic map $g: \Sigma \to S^2$.  Then, it follows from 
Fischer-Colbrie's arguments that  the index of $M$ coincides with the 
index of the Schr\"{o}dinger operator on $\Sigma$ defined by 
\[ {\cal L} = \triangle - |dg|^2 \; . \]  Here, 
$\triangle$ and $|dg|$ are 
computed relatively to any metric on $\Sigma$ that is conformally 
equivalent to the original metric on $M$.  

The purpose of this paper is to extend this circle of ideas to 
constant mean curvature 1 surfaces in $\bfH^3(-1)$.  In this case, the 
role played by 
the map $g$ is replaced by the so-called secondary Gauss map $G$ (we 
describe $G$ 
below; it is a certain multivalued map that comes from the 
Bryant-Weierstrass representation).  It so 
happens that if $M \in \bfH^3(-1)$ has constant mean curvature 1 and 
finite total curvature, then $M$ 
is also conformally finite, but it is {\em no} longer true, in 
general, that $G$ 
extends meromorphically across the ends.  This means, as we shall 
see, that the analysis necessary for studying the index in the 
$\bfH^3(-1)$ case is much more involved than in the $\bfR^3$ case.  
More precisely, let $ds^2$ and $K$ denote the induced metric and the 
Gaussian curvature in both cases.  A common feature here is that 
$K \leq 0$ and vanishes only at isolated points (unless $M$ is either 
a plane or a horosphere) and that $d\bar{s}^2 = 
-K ds^2$ is a spherical pseudo-metric on $M$ with conical singularities 
at these points.  

A crucial point here is to determine the behavior 
of $d\bar{s}^2$ at an end of $M$.  To this effect, 
let $z=(x,y)$ be a conformal parameter around some end so that the end 
corresponds to $z=0$.  Note that, since $M$ is complete, 
$ds^2$ certainly becomes infinite when $z \to 0$.  Also, 
since the surface 
has finite total curvature, the limiting value of $-K$ as $z \to 0$ is 
zero.  So at first sight, it is unclear what the behavior of 
$d\bar{s}^2 = -K ds^2$ is at the ends.  
It is well known, however, that in the minimal case we 
can choose $z$ so that 
\[ d\bar{s}^2 = \frac{4|g^\prime|^2}{(1+|g|^2)^2} |dz|^2 \; , \] 
where primes denote derivative with respect to $z$.  
Moreover, since $g$ extends meromorphically across the ends, we can 
assume $g \approx z^\ell, \ell \in \bfZ^+$, so that 
$d\bar{s}^2$ is bounded 
around the end.  One can now take advantage of this fact when one 
does the analysis necessary for examining the index of the minimal 
case, eventually obtaining 
Fischer-Colbrie's results.  In the hyperbolic case, 
we shall compute below that 
\[ d\bar{s}^2 = \frac{4|G^\prime|^2}{(1+|G|^2)^2} |dz|^2 \; . \]  But 
now we can only assume that $G(z) \approx z^\mu$, for some 
$\mu > 0$ depending only upon the end.  In 
particular, if $0<\mu<1$ for some end, $d\bar{s}^2$ is {\em not} 
bounded at this end and the analysis for the minimal case does not 
apply to this situation.  

Doing the necessary extra analysis is the heart of this paper.  
More precisely, 
we show that the canonical form for $d\bar{s}^2$ 
(just above) around the ends implies 
that the Sobolev space $\bar{H}^1 = H_{d\bar{s}^2}^1$ is {\em compactly} embedded 
in $L^2_{d\bar{s}^2}$ (Lemma 4.4).  Once this has been established, 
it is an easy matter to use standard variational methods to define 
Ind$(\Sigma)$ as being the index of a certain operator $\bar{L}$ 
defined on 
$\Sigma$ and corresponding to the Schr\"{o}dinger operator $\cal L$ 
in the 
minimal case (section 5).  It follows easily from the construction that Ind$_u 
(M) \leq $ Ind$(\Sigma )$, where Ind$_u (M)$ denotes the {\em 
unconstrained} index of $M$, namely, the index as computed for not 
necessarily volume preserving variations (it follows from our arguments 
that Ind$(M)$ and Ind$_u(M)$ differ at most by 1, so that computing 
Ind$_u(M)$ takes us a great deal of the way 
toward computing Ind$(M)$, our ultimate 
concern here).    Furthermore, using standard results in elliptic 
regularity theory, we show that the eigenfunctions of $\bar{L}$ extend 
continuously 
across the ends (Lemma 5.2).  This extra regularity property enables us to show 
that Ind$(\Sigma) \leq $ Ind$_u (M)$, after an argument due to 
Fischer-Colbrie (Lemma 5.3).

Once this analysis is done, we find that we have an alternate proof 
of the result by do Carmo and Silveira \cite{CS} 
that if a constant mean curvature 
1 surface in $\bfH^3(-1)$ has finite total curvature, then it has 
finite index (Corollary 5.2).  The advantage of our way of proving this result is that 
it gives us tools that allow us to compute explicit bounds on index 
for some concrete examples.  For some surfaces we can even 
compute the index exactly.  

For example, using methods similar to 
those of Nayatani \cite{N1}, we can compute the index of the catenoid 
cousins and the Enneper cousins of higher winding order, as well as 
some other examples described by Umehara and Yamada \cite{UY1} (section 6).  Our 
results about the index of these examples yields some surprising 
differences from the index of minimal surfaces in $\bfR^3$.  For 
example, unlike the minimal catenoid in $\bfR^3$, the catenoid 
cousins in $\bfH^3(-1)$ can have arbitrarily high index (Theorem 6.1).  Also, 
although the only minimal surfaces in $\bfR^3$ with index 1 are the 
catenoid and Enneper's surface, there are many more examples of index 
1 surfaces in the hyperbolic case (final remark of section 6).  

As another example, using methods similar to those of Choe \cite{Cho}, 
we can compute lower bounds for many constant mean curvature 1 
surfaces in $\bfH^3(-1)$ (Theorem 7.1).  We find lower bounds for genus 
1 $n$-noid cousins (Corollary 7.3), and for 
genus $k$ Costa surface cousins (Corollary 7.1).  And, in 
general, for those constant mean curvature 1 surfaces that are 
constructed via a deformation method \cite{RUY} from minimal surfaces 
in $\bfR^3$, we can find a lower bound for index (Theorem 8.1).  

The second author owes special thanks to Shin Nayatani for many 
helpful discussions.  Thanks are also due to Pierre Berard, 
Etienne Sandier, Shin Kato, and David Goldstein.  

\section{Definition of index}

Let $\Phi: M \to M^3(a)$ be an isometric 
immersion of a 2-dimensional manifold $M$ into a complete 
simply-connected 3-dimensional manifold $M^3(a)$ with constant sectional 
curvature $a$.  Let 
$\vec{N}$ be a unit normal vector field on 
$\Phi (M)$ (we write $\Phi^* \vec{N}$ simply as $\vec{N}$ defined on 
$M$).  Let $\Phi (t)$ be a smooth variation of immersions 
for $t \in (-\epsilon,\epsilon)$ so that $\Phi (0) = \Phi$.  
Assume that the variation has compact support.  We can assume that 
the corresponding variation vector field at time $t=0$ is $u
\vec{N}$, $u \in C_0^\infty(M)$.  Let $A(t)$ be the area of
$\Phi (t) (M)$ and $H$ be the mean curvature of $\Phi (M)$.  
The first variational formula (\cite{L}) is 
\[ \left. \frac{dA}{dt} \right|_{t=0} = - \int_M \langle  n H 
\vec{N}, u \vec{N} \rangle dA \; \; , \] where $\langle , \rangle$ 
and $dA$ are the metric and area form on $M$ induced by the 
immersion $\Phi$.  If $H$ is constant, then 
$A^\prime(0) = -nH \int_M
u dA$.  Let $V(t)$ be the volume of $\Phi (t) (M)$, then $V^\prime(0) = 
\int_M u dA$.  A variation is said to be {\em volume preserving} if 
$\int_M u dA = 0$.  It follows that
$\Phi (M)$ is critical for area amongst all volume preserving variations.

The second variation formula for volume preserving variations (\cite{Che}, 
\cite{Si}, \cite{L}) is 
\[ \left. \frac{d^2A}{dt^2} \right|_{t=0} := 
\int_M \{|\nabla u|^2 - 2(2a+2H^2-K)u^2\} dA \; , \] where $K$ is the 
Gaussian curvature on $M$.  
Since we will be investigating surfaces of constant mean curvature 1 in 
hyperbolic space with constant sectional curvature $-1$, we will 
restrict ourselves to the case $a = -1$ and $H = 1$, so 
\[ \left. \frac{d^2A}{dt^2} \right|_{t=0} = 
\int_M \{|\nabla u|^2 + 2 K u^2\} dA \; . \]
This formula is the same for both minimal surfaces 
in $\bfR^3:=M^3(0)$ and 
constant mean curvature 1 surfaces in $\bfH^3:=M^3(-1)$, giving us our 
first indication of the close relationship between these two types of 
surfaces.  Another indication of this close relationship is the 
Weierstrass representations described in the next section.  

The index Ind($M$) is the maximum possible 
dimension of a subspace of volume preserving variation functions 
in $C_0^\infty(M)$ on which $\left. \frac{d^2A}{dt^2} \right|_{t=0} < 0$.  
The purpose of this paper is to estimate Ind($M$).  

We define Ind$_u$($M$) as the maximum possible 
dimension of a subspace of (not necessarily volume preserving) 
variation functions in $C_0^\infty(M)$ on which the above 
$\left. \frac{d^2A}{dt^2} \right|_{t=0} < 0$.  
(The subscript $u$ stands for ``{\em u}nconstrained index''.)  
Clearly, Ind$_u$($M$) $\geq$ Ind($M$).  
We will show later that 
also Ind$_u$($M$) $- 1 \leq$ Ind($M$).  
The methods we use in this 
paper allow us to compute Ind$_u$($M$), but what we really want to 
compute is Ind($M$).  However, these two indices can differ by 
at most 1, so computing Ind$_u$($M$) means that we know 
Ind($M$) must be either Ind$_u$($M$) 
or Ind$_u$($M$)$-1$.  

\section{The Weierstrass representation}

Both minimal surfaces in $\bfR^3$ and constant mean curvature 
1 surfaces in 
$\bfH^3$ can
be described parametrically by a pair of meromorphic functions on a 
Riemann surface, via a Weierstrass representation.  First we describe 
the well-known Weierstrass representation for minimal surfaces in $\bfR^3$.
We will incorporate into this representation 
the fact that any complete minimal surface of 
finite total curvature is 
conformally equivalent to a Riemann surface 
$\Sigma$ with a finite number of 
points $\{p_j\}_{j=1}^k \subset 
\Sigma$ removed (\cite{O}):  

\begin{lemma}
Let $\Sigma$ be a Riemann surface.  Let $\{p_j\}_{j=1}^k \subset 
\Sigma$ be a finite number of points, which will represent the ends 
of the minimal surface defined in this lemma.  
Let $z_0$ be a fixed point in $\Sigma \setminus \{p_j\}$.
Let $g$ be a meromorphic function from $\Sigma \setminus \{p_j\}$ to 
$\bfC$.  
Let $f$ be a holomorphic function from $\Sigma \setminus \{p_j\}$ 
to $\bfC$.  Assume 
that, for any point in $\Sigma \setminus \{p_j\}$, $f$ has a zero of 
order $2k$ at some point if and only if $g$ has 
a pole of order $k$ at that point, and assume that $f$ has no other 
zeroes on $\Sigma \setminus \{p_j\}$.  Then 
\[ \Phi(z) = \mbox{Re} \int_{z_0}^{z}
\; \left( \begin{array}{c}
        (1-g^2)f d\zeta \\
        i(1+g^2)f d\zeta \\
        2gf d\zeta
        \end{array}
\right) \] is a conformal minimal immersion of the universal cover 
$\widetilde{\Sigma \setminus \{p_j\}}$ of 
$\Sigma \setminus \{p_j\}$ into $\bfR^3$. 
Furthermore, any complete minimal surface with finite total 
curvature in $\bfR^3$ can be represented in this way.
\end{lemma}

The map $g$ can be geometrically interpreted as the stereographic 
projection of the Gauss map.  
The first and second fundamental forms and the intrinsic Gaussian 
curvature for the surface $\Phi$ are 
    \[ ds^2 = (1+g\bar g)^2\,fdz \cdot \overline{fdz} \; , \; \; 
        II = -2\mbox{Re}(Q) \; , \; \; 
   K = -4 \left( \frac{|g^\prime|}{|f| (1+|g|^2)^2} \right)^2  \; , \]
where the Hopf differential $Q$ is defined to be $Q = g^\prime f dz^2$.  

To make a surface of finite total curvature (i.e. $\int_\Sigma -K dA 
< + \infty$, which is necessary 
to make a surface of finite index \cite{FC}) we must choose $f$ and
$g$ so that $\Phi$ is well defined on $\Sigma 
\setminus \{p_j\}$ itself.  Usually this 
involves adjusting some real parameters in the descriptions of $f$ 
and $g$ and $\Sigma \setminus \{p_j\}$ so that the 
real part of the above integral about any
nontrivial loop in $\Sigma \setminus \{p_j\}$ is zero.

We now describe a Weierstrass type representation for 
constant mean curvature $c$
surfaces in $\bfH^3(-c^2):=M^3(-c^2)$.  
This result is a composite of several
results that are found in \cite{By}, \cite{UY3}, \cite{UY4}.  

\begin{lemma}
Let $\Sigma$, $\Sigma \setminus \{p_j\}$, 
$z_0$, $f$, and $g$ be the same as in the previous 
lemma.  Choose a null holomorphic immersion $F:
\widetilde{\Sigma \setminus \{p_j\}} \to SL(2,\bfC)$ so 
that $F(z_0)$ is the identity
matrix and so that $F$ satisfies 
\begin{equation}
F^{-1}dF = c \left( \begin{array}{cc} 
g & -g^2 \\ 
1 & -g 
\end{array} \right) f dz \; , 
    \label{eq:wode}
\end{equation}
then 
$\Phi:\widetilde{\Sigma \setminus \{p_j\}} \to H^3(-c^2)$ defined by 
\begin{equation}
        \Phi = \frac{1}{c} F^{-1} \overline{F^{-1}}^t
    \label{eq:imm}
\end{equation}
is a conformal constant mean curvature $c$ 
immersion into $\bfH^3(-c^2)$ with the Hermitean model.
Furthermore, any constant mean curvature $c$ surface with finite 
total curvature in $\bfH^3(-c^2)$ can be represented in this way.
\end{lemma}

We call $g$ the {\it hyperbolic Gauss map} of $\Phi$.   
As its name suggests, 
the map $g(z)$ has a geometric interpretation for this case as well.  
It is the image
of the composition of two maps.  The first map is from each point on 
the surface to the point at the sphere at infinity in the Poincare 
model
which is at the opposite 
end of the oriented perpendicular geodesic ray starting at the point
$z$ on the surface.  The second map is stereographic projection of the
sphere at infinity to the complex plane $\bfC$ \cite{By}.
The first and second fundamental 
forms and the intrinsic Gaussian curvature 
of the surface are 
    \[
        ds^2 = (1+G\bar G)^2 \, \frac{f g^\prime}{G^\prime} 
         \overline{\left(\frac{f g^\prime}{G^\prime}\right)} dz 
         \overline{dz} \; , \; \; 
        II = -2 \mbox{Re}(Q) + c\,ds^2 \; , \; \; 
   K = -4 \left( \frac{|G^\prime|^2}{|g^\prime| 
|f| (1+|G|^2)^2} \right)^2  \; ,
    \label{eq:second}
    \] where in this case the Hopf differential is $Q = - f g^\prime 
dz^2$ (the sign change in $Q$ is due to the fact that we are 
considering the ``dual'' surface; see \cite{UY4} for an 
explanation of this), and 
where $G$ is defined as the multi-valued meromorphic function 
\[
G=\frac {dF_{11}}{dF_{21}}=\frac {dF_{12}}{dF_{22}} 
\]
on $\Sigma \setminus \{p_j\}$, 
with $F=(F_{ij})_{i,j=1,2}$.  The reason that $G$ is multi-valued is 
that $F$ itself can be multi-valued on $\Sigma \setminus \{p_j\}$ 
(even if $\Phi$ is well defined on $\Sigma \setminus \{p_j\}$ itself).
The function $G$ is called the {\it secondary Gauss map}
of $\Phi$ (\cite{By}).  

In the above lemma, we have 
changed the notation slightly from the notation used 
in $\cite{By}$ and $\cite{RUY}$, 
because we wish to use the same symbol ``$g$'' both for the map $g$ 
used in the Weierstrass representation for minimal surfaces in 
$\bfR^3$ and for the hyperbolic Gauss map used in the Weierstrass 
representation for constant mean curvature surfaces in 
$\bfH^3$.  And we further wish to give a separate notation ''$G$'' for 
the secondary Gauss map used in the hyperbolic case.  We do this to 
emphasize that, in relation to their geometric interpretations, the 
``$g$'' in the 
Euclidean case is more closely related to the hyperbolic Gauss map 
``$g$'' in the $\bfH^3$ case than 
to the secondary Gauss map ``$G$'' (as we will see in section 6).  

In order for $\Phi$ to be well-defined on $\Sigma \setminus \{p_j\}$ 
itself, 
it is sufficient 
and necessary that $F$ satisfy a condition called the 
$SU(2)$-{\it condition}.  Note that if one travels about a nontrivial 
loop in $\Sigma \setminus \{p_j\}$, 
then $F \to B F$, where $B \in SL(2,\bfC)$.  If for every
loop in $\Sigma \setminus \{p_j\}$, the resulting matrix 
$B$ satisfies $B \in SU(2)$, then 
the
$SU(2)$-condition is satisfied.  If $B \in SU(2)$, then 
$F^{-1} \overline{F^{-1}}^t = 
(B F)^{-1} \overline{(B F)^{-1}}^t$, so it follows that if 
the 
$SU(2)$-condition holds, then $\Phi$ is well defined on 
$\Sigma \setminus \{p_j\}$ itself.
When $F \to B F$, we have the following effect on the secondary 
Gauss map:  
\[ G \rightarrow \frac{b_{11} G + b_{12}}{b_{21} G + b_{22}} \; , 
\mbox{ for } B = \left(
\begin{array}{cc} b_{11} & b_{12} \\ b_{21} & b_{22} \end{array}
\right) \in SU(2) \; . \]

We now state some known facts, which when taken together, show that 
constant mean curvature 1 surfaces in $\bfH^3$ and minimal surfaces 
in $\bfR^3$ are very closely related.  These facts provide the 
motivation for the results in sections 7 and 8 of this paper:

\begin{itemize}
\item It was shown in \cite{UY2} that if $f$ and $g$ and 
$\Sigma \setminus \{p_j\}$ are fixed, 
then as $c \to 0$, the constant mean curvature $c$ 
surfaces $\Phi$ in $\bfH^3(-c^2)$ 
converge to a minimal surface in $\bfR^3$.  This can be sensed 
from the fact that $G \to g$ and $B \to$ identity as $c \to 0$ (which 
follow directly from equation \ref{eq:wode}), and hence 
the above first 
and second fundamental forms for the constant mean curvature $c$ 
surfaces $\Phi$ converge to the fundamental forms for a 
minimal surface as $c \to 0$ (up to a sign change in $II$ -- a 
change of orientation).  
\item It was shown in \cite{RUY} that a finite total curvature minimal
surface in $\bfR^3$ satisfying certain nondegeneracy and symmetry 
conditions (these conditions 
are fairly general and include most known examples) 
can be deformed into a constant mean curvature $c$ surface in 
$\bfH^3(-c^2)$ for $c \approx 0$, 
so that $\Sigma$, $f$, and $g$ are the same, up to 
a slight adjustment of the real parameters that are used to solve the 
period problem.  The deformed surface might not have finite total 
curvature, but it will be of the same topological type as the minimal 
surface, and it will have the same reflectional symmetries as the 
minimal surface.  
\item Consider the Poincare model for $\bfH^3(-c^2)$ for $c \approx 0$.  
It is a round ball in $\bfR^3$ centered at the origin 
with Euclidean radius $\frac{1}{c}$ endowed with a complete 
radially-symmetric metric 
$ds_c^2 = \frac{4 \sum dx_i^2}{(1-c^2 \sum x_i^2)^2}$ 
of constant sectional curvature $-c^2$.  Contracting 
this model by a factor of $c$, we obtain a map to the Poincare model 
for $\bfH^3$.  Under this mapping, constant mean curvature $c$ 
surfaces are mapped to
constant mean curvature $1$ surfaces.  Thus the problem of existence of 
constant mean curvature $c$ surfaces in $\bfH^3(-c^2)$ 
for $c \approx 0$ is equivalent to the problem of existence of 
constant mean curvature $1$ surfaces in $\bfH^3$.  Furthermore, under 
this mapping, the area form on the constant mean curvature $c$ surface 
is changed only by a constant 
factor $c^2$: If $dA_c$ is the area form on the constant mean 
curvature $c$ surface, and $dA_1$ is the area form on the constant mean 
curvature $1$ surface, then $dA_c = c^2 dA_1$.  
Hence a variation that reduces area on the constant 
mean curvature $c$ surface is mapped to a variation that 
reduces area on the constant mean curvature $1$ surface (and 
vice-versa).  Hence this mapping preserves the index.  
\end{itemize} 

\section{Showing that $\bar{H}^{1}$ is compactly contained in 
$L^{2}_{d\bar{s}^2}$}

We now consider $M$ to be a complete 
constant mean curvature 1 surface in 
$\bfH^3$ with finite total curvature.  We will assume the surface is 
not a horosphere.  (Assuming that the surface is not a horosphere will 
not add any extra conditions to our index results, since the index of 
the horosphere is known to be zero \cite{Si}.)  Suppose that 
$M$ has Weierstrass representation $\Phi: \Sigma \setminus \{p_j\} 
\to M$ 
with Riemann surface $\Sigma$ and functions $f,g:\Sigma \to \bfC$, and 
that $G$ is the secondary Gauss map.  Let $ds^2$ be the complete metric 
on $M$ pulled back to $\Sigma$.  Note that $M$ is conformally 
equivalent to $\Sigma$ 
with a finite number of points $\{p_j\}$ removed; each removed point 
$p_j$ corresponds to an end of $M$.  So $ds^2$ is 
defined on $\Sigma \setminus \{p_j\}$.  Let $d\bar{s}^2 = 
G^* ds^2_{S^2} = -K ds^2$ be the 
singular pull back metric of the canonical metric on $S^2$ via the 
secondary Gauss map $G$, defined on $\Sigma \setminus \{ p_j \}$, but 
with isolated singularities where $K = 0$.  We let 
$d\tilde{s}^2$ be a conformal nonsingular metric defined on 
$\Sigma$.  Any choice for $d\tilde{s}^2$ will suffice, provided it 
is conformally equivalent to $ds^2$ on $\Sigma \setminus \{p_j\}$.  
Let $dA$ (resp. $d\bar{A}$, $d\tilde{A}$) and 
$\nabla$ (resp. $\bar{\nabla}$, $\tilde{\nabla}$) and 
$\triangle$ (resp. $\bar{\triangle}$, $\tilde{\triangle}$) be the 
area form and gradient and Laplacian on $\Sigma$ 
with respect to the metric $ds^2$ (resp. 
$d\bar{s}^2$, $d\tilde{s}^2$).  

We choose the sign of the Laplacian so that $\int_\Omega 
|\nabla u|^2 = + \int_\Omega u \triangle u$ for any $u \in 
C^\infty_0(\Omega)$.  
(Thus, for example, the Laplacian on the standard Euclidean plane 
$\bfR^2$ will be 
$-\frac{\partial^2}{\partial x^2}-\frac{\partial^2}{\partial y^2}$.)  
So if $u \in C^\infty_0(\Omega)$ satisfies $\triangle u = \lambda u$ 
for some constant $\lambda$, 
where $\Omega$ is a region in $\Sigma$, then $\int_\Omega 
|\nabla u|^2 = \int_\Omega u \triangle u = \lambda \int_\Omega u^2$ 
and so $\lambda \geq 0$.  Thus our convention for the sign of the 
Laplacian implies that the 
eigenvalues of the Laplacian will be nonnegative.  

We now list some easily determined facts that will be used 
throughout this and the next section.  We can define 
$|dG|_{ds^2}^2$ (resp. $|dG|_{d\bar{s}^2}^2$, 
$|dG|_{d\tilde{s}^2}^2$) by 
$|dG|_{ds^2}^2 = \sum_{j=1}^2 \left< dG(e_j), dG(e_j) \right>_{ds^2_{S^2}}$ 
(resp. $|dG|_{d\bar{s}^2}^2 = \sum_{j=1}^2 \left< dG(\bar{e}_j), 
dG(\bar{e}_j) \right>_{ds^2_{S^2}}$, 
$|dG|_{d\tilde{s}^2}^2 = \sum_{j=1}^2 \left< dG(\tilde{e}_j), 
dG(\tilde{e}_j) \right>_{ds^2_{S^2}}$), 
where $dG$ is the tangent map of $G$ and $\{e_1,e_2\}$ 
(resp. $\{\bar{e}_1,\bar{e}_2\}$, $\{\tilde{e}_1,\tilde{e}_2\}$) 
is an orthonormal 
basis of vector fields with respect to the metric 
$ds^2$ (resp. $d\bar{s}^2$, $d\tilde{s}^2$).  The following hold:  
\begin{itemize}
\item $d\bar{s}^2 = -K ds^2$, $d\bar{A} = -K dA$,   
$\triangle = - K \bar{\triangle}$
\item $d\bar{s}^2 = \frac{1}{2} |dG|_{d\bar{s}^2}^2 d\bar{s}^2 
= \frac{1}{2} |dG|_{ds^2}^2 ds^2
= \frac{1}{2} |dG|_{d\tilde{s}^2}^2 d\tilde{s}^2\; \; \; \; $ 
(conformal invariance).  
\end{itemize}

We now consider the variation described in the second section with 
variation vector field $u \vec{N}$ on $M$ at time $t=0$.  Since 
$\left. \frac{d^2A}{dt^2} \right|_{t=0}$ obviously depends on $u$, we will 
write it as $\left. \frac{d^2A}{dt^2} \right|_{t=0}(u)$.  In the next lemma, 
we will consider $\Sigma$ and $u$ to be 
fixed, but we consider whether or not 
$\left. \frac{d^2A}{dt^2} \right|_{t=0}(u)$ depends on $G$, $g$, and $f$.

We now state a crucial computation -- it is crucial because it 
explains why the pull-back of the metric on the sphere via the map 
$G$ plays such a dominant role in computing Ind($M$), and explains why 
the operators $L$ and $\bar{L}$ (defined later) are somehow 
``the same'' operator: 
\[ \left. \frac{d^2A}{dt^2} \right|_{t=0}(u) = \int_\Sigma 
\{|\nabla u|^2 + 2 K u^2\} dA = \int_\Sigma 
\{u\triangle u + 2Ku^2\} dA = \]\[
\int_\Sigma \{ -uK\bar{\triangle} u + 2Ku^2\} dA = 
\int_\Sigma \{ -u \bar{\triangle} u + 2u^2\} KdA = 
\int_\Sigma \{ u \bar{\triangle} u - 2u^2\} d\bar{A} \; . \] 
Since the integrand $\{ u \bar{\triangle} u - 2u^2\} d\bar{A}$ 
is completely determined by the pull-back of the spherical metric via 
the map $G$, we know that $\left. 
\frac{d^2A}{dt^2} \right|_{t=0}(u)$ depends 
only on $G$, and does not depend on $g$ and $f$.  

\begin{lemma}
$\left. \frac{d^2A}{dt^2} \right|_{t=0}(u)$ is completely independent of 
$f$ and $g$.  It does depend on $G$, but not on the choice of value of the 
multi-valued $G$.  
\end{lemma}

\begin{proof}
As noted above, $\left. 
\frac{d^2A}{dt^2} \right|_{t=0}(u)$ depends 
only on $G$, not on $g$ and $f$.  
Clearly, $\left. \frac{d^2A}{dt^2} \right|_{t=0}(u)$ does depend on $G$, but 
to show that it does not depend on the choice of value of the 
multi-valued $G$, we first show that the 
first fundamental form is independent of the $SU(2)$-condition.
Let $\hat{f} = -g^\prime f/G^\prime$.  Travelling about a 
loop in $\Sigma$ corresponding to a homologically nontrivial loop in 
$M$, we have $F \to BF$ for some $B \in $ SU(2), and as we saw before 
$G \rightarrow (b_{11} G + b_{12})/(b_{21} G + b_{22})$, where 
$b_{ij}$ are the entries of $B$.  Thus travelling about the loop 
makes the transformation $G^\prime \to 
G^\prime/(b_{21} G + b_{22})^2$, 
and since $f$ and $g$ are left unchanged (and 
therefore $g^\prime$ is also unchanged), it follows that 
$\hat{f} \to \hat{f} (b_{21} G + b_{22})^2$.
Now we consider the conformal factor 
$\hat{f} \overline{\hat{f}} (1+G \bar{G})^2)^2$ 
in the first fundamental form.  Denoting $(b_{11} G + 
b_{12})/(b_{21} G + b_{22})$ as $B \cdot G$, we see that 
\[ \hat{f} \overline{\hat{f}} (1+G \bar{G})^2)^2 \to
(b_{21} G + b_{22})^2 \hat{f} \overline{(b_{21} G + b_{22})^2 \hat{f}}
(1 + (B \cdot G) (\overline{B \cdot G}))^2 
= \hat{f} \overline{\hat{f}} (1+G \bar{G})^2)^2 \; , \]
since $b_{11} = \overline{b_{22}}$ 
and $b_{12} = - \overline{b_{21}}$.

Therefore $\langle , \rangle$ is independent of the 
SU(2)-condition, and therefore $\nabla u$ and $dA$ are independent
of the $SU(2)$-condition, since they are determined by the first 
fundamental form.  And since $K$ depends only on the first
fundamental form, $K$ is also independent of the $SU(2)$-condition.
We conclude that $\left. \frac{d^2A}{dt^2} \right|_{t=0}(u)$ is 
independent of the $SU(2)$-condition.  Thus 
$\left. \frac{d^2A}{dt^2} \right|_{t=0}(u)$ is well-defined even 
though $G$ is multi-valued.
\end{proof}

For any $p \geq 2$, let $L^p_{d\tilde{s}^2}(\Omega)$ (resp. 
$L^p_{d\bar{s}^2}(\Omega)$) be the space of 
measurable functions $f$ on 
$\Omega \subset \Sigma$ such that $\int_{\Omega} 
|f(x)|_{d\tilde{s}^2}^p d\tilde{A} < \infty$ (resp. 
$\int_{\Omega} 
|f(x)|_{d\bar{s}^2}^p d\bar{A} < \infty$).  
In the case that $\Omega = \Sigma$, we may write simply 
$L^p_{d\tilde{s}^2}$ (resp. 
$L^p_{d\bar{s}^2}$) instead of 
$L^p_{d\tilde{s}^2}(\Sigma)$ (resp. 
$L^p_{d\bar{s}^2}(\Sigma)$).  

We now begin to work toward a proof that $\bar{H}^{1}$ is compactly 
contained in $L^{2}_{d\bar{s}^2}$.  

\begin{lemma}
If $p$ is sufficiently large, then 
$L^p_{d\tilde{s}^2}$ is continuously 
contained in both $L^2_{d\tilde{s}^2}$ and 
$L^2_{d\bar{s}^2}$.  
\end{lemma}

\begin{proof}
Since $\Sigma$ is compact, $L^p_{d\tilde{s}^2}$ is 
continuously contained in $L^2_{d\tilde{s}^2}$ for all $p \geq 2$.  
(See, for example, \cite{GT}, equation (7.8)).  

As for the second assertion, consider a point $p_j \in \Sigma$ 
representing an end of the complete surface.  Let 
$U_j$ be a small neighborhood of $p_j$.  We may choose $d\tilde{s}^2$ 
so that $d\tilde{s}^2 = dx^{2} + dy^{2} = 4dzd\bar{z}$ on $U_j$.  
We now show that locally on $U_j$, 
\[ d\bar{s}^2 \approx 4 \mu^2 \frac{r^{2 \mu - 2}}
{(1+r^{2 \mu})^2} d\tilde{s}^2 \; , \] 
with $r = \sqrt{x^2+y^2}$.  (The symbol "$\approx$" means that for 
functions $a(z)$, $b(z)$ defined in a neighborhood of $z=0$, 
$a(z) \approx b(z)$ if for all $\epsilon > 0$, there exists a 
$\delta > 0$ such that $|z| < \delta$ implies $|
\frac{a(z)}{b(z)} - 1| < \epsilon$.)  
The above relation follows from the fact 
that locally near an end we can make the following 
normalization: we can choose the complex coordinate $z$ on $U_j$ so 
that the end $p_j$ is at $z=0$.  
By the previous lemma, we may change $G$ to 
$(b_{11} G + b_{12})/(b_{21} G + b_{22})$ for any 
$B = \{b_{ij}\} \in $ SU(2), without affecting the second variation 
formula.  We may choose $B$ so that 
$((b_{11} G + b_{12})/(b_{21} G + b_{22}))(z=0) = 0$.  Hence we may 
assume that $G(0) = 0$.  We then have that $G = 
z^\mu \hat{G}$, where $\hat{G}$ is a holomorphic function in a 
neighborhood of $z=0$ such that $\hat{G}(0) \neq 0$, for some $\mu \in 
\bfR^+$, where $z = x+iy$ \cite{UY1}.  Changing $z$ to 
$\hat{G}(0)^{-1/\mu} z$ if necessary, we may assume that $\hat{G}(0) = 
1$.  

The point 
corresponding to $G(z)$ under the inverse of stereographic projection 
is \[ {\cal G} = ({\cal G}_1,{\cal G}_2,{\cal G}_3) = 
\frac{1}{|G|^{2}+1}(2\mbox{Re}(G),2\mbox{Im}(G),|G|^2 - 1) \; . \]
Note that for any real-valued function $f:\bfC \rightarrow \bfR$, 
we have $f_{z} = \frac{1}{2} (f_{x} - i f_{y})$, so $|f_{z}|^{2} = 
\frac{1}{4} (|f_{x}|^{2} + |f_{y}|^{2})$.  For any complex-valued 
holomorphic function $f:\bfC \rightarrow 
\bfC$, we have $\frac{\partial}{\partial z} (\mbox{Re}(f)) = 
\frac{\partial}{\partial z} \frac{1}{2} (f + \bar{f}) = 
\frac{1}{2} f_{z}$, 
$\frac{\partial}{\partial z} (\mbox{Im}(f)) = 
\frac{\partial}{\partial z} \frac{i}{2} (\bar{f} - f) = 
\frac{-i}{2} f_{z}$, and 
$\frac{\partial}{\partial z} (f \bar{f}) = 
f_{z} \bar{f} + f \bar{f}_{z} = \bar{f} f_{z}$.  Using these 
properties and the fact that $G = z^\mu \hat{G}$, we have 
\[ |dG|^{2}_{d\tilde{s}^2} = 
|d{\cal G}|^{2}_{d\tilde{s}^2} = 
\sum_{i=1}^{3} |\nabla {\cal G}_{i}|^{2} = 
\sum_{i=1}^{3} (|({\cal G}_{i})_{x}|^{2} + |({\cal G}_{i})_{y}|^{2}) 
\]\[ = 
\sum_{i=1}^{3} 4 |({\cal G}_{i})_{z}|^{2} = \frac{8G_{z}\overline{G_z}}
{(1+G\bar{G})^{2}} \approx \frac{8 
\mu^{2}r^{2\mu-2}}{(1+r^{2\mu})^{2}} \; . \]  
And thus it follows that $d\bar{s}^2 \approx 
4 \mu^2 \frac{r^{2 \mu - 2}}
{(1+r^{2 \mu})^2} d\tilde{s}^2$ on $U_j$.  

Suppose $u \in L^p_{d\tilde{s}^2}(U_j)$.  By the 
Holder inequality we have $\int_{U_j} u^2 d\bar{A} = 
\int_{U_j} u^2 \frac{1}{2} |dG|^{2}_{d\tilde{s}^2} d\tilde{A} 
\leq A^{\frac{2}{p}} 
B^{\frac{1}{q}}$, where $A = \int_{U_j} u^{p} 
d\tilde{A}$ and 
$B < c \int_{U_j} r^{(2 \mu - 2) q} d\tilde{A}$, with 
$\frac{2}{p} + \frac{1}{q} = 1$ and $c>0$ some finite constant.  
If $q$ is close enough to 1, 
then $B$ is finite, since $\mu > 0$ and $d\tilde{A}$ has the local 
expression $d\tilde{A} = rdrd\theta$ in polar coordinates.  So 
there exists a constant $k_j$ such that $||u||_{L^2_{d\bar{s}^2}(U_j)} 
\leq k_j ||u||_{L^p_{d\tilde{s}^2}(U_j)}$ for each $j$.  

On $\Sigma \setminus \{\cup U_j\}$, $d\bar{s}^2$ is bounded.  So it 
is clear from the Holder inequality that there exists a constant $k_0$ 
such that $||u||_{L^2_{d\bar{s}^2}(\Sigma \setminus \{\cup U_j\})} 
\leq k_0 ||u||_{L^p_{d\tilde{s}^2}(\Sigma \setminus \{\cup U_j\})}$.  
Let $k = \max\{k_0,k_j\}$.  Then, choosing $p$ large enough, we have 
$||u||_{L^2_{d\bar{s}^2}} 
\leq k ||u||_{L^p_{d\tilde{s}^2}}$.  
\end{proof}

\begin{remark}
If $\mu \geq 1$ for all ends, then $d\bar{s}^2$ is bounded on all of 
$\Sigma$, and the lemma 
holds even for $p = 2$.  We could argue this way: 
suppose $u \in L^{2}_{d\tilde{s}^2}$.  Then 
$\int_{\Sigma} u^2 d\bar{A} = 
\int_{\Sigma} u^2 \frac{1}{2} 
|dG|_{d\tilde{s}^2}^2 d\tilde{A} \leq
\mbox{(const)} \int_{\Sigma} u^{2} d\tilde{A}$.  
Thus $||u||_{L^{2}_{d\bar{s}^2}} \leq 
\mbox{(const)} ||u||_{L^{2}_{d\tilde{s}^2}}$.  So 
$L^{2}_{d\tilde{s}^2}$ is 
continuously included in $L^{2}_{d\bar{s}^2}$.  
\end{remark}

We define $\tilde{H}^{1}(\Sigma) = 
\{ u \in L^2_{d\tilde{s}^2} \; | \; 
du \in L^2_{d\tilde{s}^2} \}$, where the derivative $du = 
(\frac{\partial u}{\partial x_1},\frac{\partial u}{\partial x_2})$ 
satisfies $\int_{\Sigma} \langle \frac{\partial u}{\partial x_{i}}, \phi 
\rangle_{d\tilde{s}^2} d\tilde{A} = \int_{\Sigma} 
\langle u, \frac{\partial \phi}{\partial x_{i}} 
\rangle_{d\tilde{s}^2} d\tilde{A}$ 
for all test functions $\phi \in C^{\infty}(\Sigma)$ and all 
coordinate functions $x_{i}$.  The condition that $du$ must satisfy 
depends on $d\tilde{s}^2$, but it is well known that $\tilde{H}^{1}$ is 
independent of $d\tilde{s}^2$ if $d\tilde{s}^2$ is 
a true metric and not a pseudometric.  We define 
$\bar{H}^{1}(\Sigma) = \{ u \in 
L^2_{d\bar{s}^2} \; | \; 
du \in L^2_{d\bar{s}^2} \}$, where $du$ satisfies $\int_{\Sigma} 
\langle \frac{\partial u}{\partial x_{i}}, \phi 
\rangle_{d\bar{s}^2} d\bar{A} = \int_{\Sigma} 
\langle u, \frac{\partial \phi}{\partial x_{i}} 
\rangle_{d\bar{s}^2} d\bar{A}$ 
for all test functions $\phi \in C^{\infty}(\Sigma)$ and all 
coordinate functions $x_{i}$.  
Note that $d\bar{s}^2$ is a psuedometric and might not be a true 
metric even away from the ends $p_j$ of the surface, since the secondary 
Gauss map may have branch points even at finite points on the surface.  
We define the two norms 
\[ \tilde{||} u \tilde{||}^{2} := \int_{\Sigma} (|\tilde{\nabla} 
u|_{d\tilde{s}^2}^{2}            
+ u^{2}) d\tilde{A} \]
\[ \bar{||} u \bar{||}^{2} := \int_{\Sigma} (|\bar{\nabla} 
u|_{d\bar{s}^2}^{2} 
+ u^{2}) d\bar{A} = \int_{\Sigma} 
(|\tilde{\nabla} u|_{d\tilde{s}^2}^{2} 
+ \frac{1}{2} |dG|_{d\tilde{s}^2}^{2} u^{2}) d\tilde{A} \; . \]

\begin{lemma}
$\bar{H}^{1}$ is continuously contained in 
$\tilde{H}^{1}$.
\end{lemma}

\begin{proof}
We need to show that there exists a $c>0$ such that
$\tilde{||} \cdot \tilde{||} \leq c \bar{||} \cdot \bar{||}$.  

By way of contradiction, suppose that such a $c$ cannot exist.  Then there 
exists a sequence $\{u_{n}\}_{n=1}^{\infty}$ of functions such that
$\tilde{||} u_n \tilde{||} = 1$ and 
$\bar{||} u_n \bar{||} < \frac{1}{n}$.  
Note the following three facts:
\begin{itemize}
\item Any bounded sequence in a Hilbert space 
has a weakly convergent subsequence (see, for example, 
\cite{GT}, p85).  In our case 
the Hilbert space is $(\tilde{H}^{1}, \tilde{||} \cdot \tilde{||})$.  
\item The inclusion of $\tilde{H}^{1}$ into 
$L^{p}_{d\tilde{s}^2}$ is compact for all $p \in [2,\infty)$ 
(See, for example, \cite{GT}, Theorem 7.22, or see \cite{Ad}.)  
\item $\tilde{||} \cdot \tilde{||}$ is lower 
semicontinuous with respect to weak 
convergence; that is, if $u_{n} \rightarrow u$ weakly, then 
$\tilde{||} u \tilde{||} \leq \liminf_{n \rightarrow \infty} 
\tilde{||} u_n \tilde{||}$.  
\end{itemize}
By the first fact, we may assume that $\{u_{n}\}$ 
converges weakly in $\tilde{H}^{1}$ to some $u \in \tilde{H}^{1}$.  
By the second fact, we may assume that $\{u_{n}\}$ 
converges strongly in $L^2_{d\tilde{s}^2}$ to some 
$v \in L^{2}_{d\tilde{s}^2}$.  Since $\tilde{H}^{1} 
\subset L^{2}_{d\tilde{s}^2}$ continuously, we have that 
$u_{n} \rightarrow u \in L^{2}_{d\tilde{s}^2}$ 
weakly.  And since the weak limit is unique, we have $u=v$.  

By the third fact, we have
\[ \int_{\Sigma} (|\tilde{\nabla} 
u|_{d\tilde{s}^2}^{2} 
+ u^{2}) d\tilde{A} \leq \liminf_{n \rightarrow \infty} 
\int_{\Sigma} (|\tilde{\nabla} 
u_{n}|_{d\tilde{s}^2}^{2} 
+ u_{n}^{2}) d\tilde{A} \; . \]  
We have strong convergence of $u_n$ to $u$ 
in $L^{2}_{d\tilde{s}^2}$, hence 
$\int_\Sigma u_{n}^{2} d\tilde{A} \to \int_\Sigma u^{2} 
d\tilde{A}$.  So we have
\[ \int_{\Sigma} |\tilde{\nabla} u|_{d\tilde{s}^2}^{2} 
d\tilde{A} \leq \liminf_{n \rightarrow \infty} 
\int_{\Sigma} |\tilde{\nabla} u_{n}|_{d\tilde{s}^2}^{2} 
d\tilde{A} \; . \]  And then since $\int_{\Sigma} 
|\tilde{\nabla} u_{n}|_{d\tilde{s}^2}^{2} 
d\tilde{A} = 
\int_{\Sigma} |\bar{\nabla} u_{n}|_{d\bar{s}^2}^{2} 
d\bar{A} < 
\frac{1}{n}$, we have $\int_{\Sigma} 
|\tilde{\nabla} u|_{d\tilde{s}^2}^{2} 
d\tilde{A} = 0$.
Therefore $u$ is constant almost everywhere.  Since $1 - \frac{1}{n}
\leq \int_{\Sigma} u_{n}^{2} d\tilde{A} \leq 1$ we have 
$\int_{\Sigma} u^{2} d\tilde{A} = 1$  (Here again we are using 
that $u_{n} \rightarrow u$ strongly in $L^{2}_{d\tilde{s}^2}$.)  
Therefore $u$ is equal to a nonzero constant almost everywhere.  

By the previous lemma, $L^{p}_{d\tilde{s}^2}$ is continuously 
included in $L^{2}_{d\bar{s}^2}$ for $p$ large enough.  
By the second fact, $\tilde{H}^{1}$ is compactly contained in 
$L^{p}_{d\tilde{s}^2}$, so it follows that 
$\tilde{H}^{1}$ is compactly contained in $L^{2}_{d\bar{s}^2}$.  
This means that any weakly convergent sequence in $\tilde{H}^{1}$ 
(which is therefore a bounded sequence in $\tilde{H}^{1}$) has
a strongly convergent subsequence in $L^{2}_{d\bar{s}^2}$.  
So, since $\int_{\Sigma} u_{n}^{2} 
d\bar{A} \leq \frac{1}{n}$, we have $\int_{\Sigma} u^{2} 
d\bar{A} = 0$.  But $u$ is a nonzero constant, so 
$0 = \int_{\Sigma} u^{2} d\bar{A} = (\mbox{const} \neq 0) \cdot 
\int_{\Sigma} \frac{1}{2} |dG|_{d\tilde{s}^2}^{2} 
d\tilde{A}$.  Therefore $|dG|_{d\tilde{s}^2} = 0$
almost everywhere, and thus $|G^\prime| = 0$.  Hence 
$K = -4 \left(|G^\prime|^2/(|g^\prime| 
|f| (1+|G|^2)^2) \right)^2 = 0$, which implies the surface 
is umbilic, and hence a horosphere.  But we assumed the surface is not 
a horosphere, so this is a contradiction.  
\end{proof}

\begin{remark}
{\em If} $\mu \geq 1$ {\em for all ends of a constant mean curvature 1 
surface} $M${\em , then} 
$\tilde{H}^{1} = \bar{H}^{1}$.  
We already know that $\bar{H}^{1}$ is continuously 
included in $\tilde{H}^{1}$, so to show this it remains only to show 
that there exists a $c>0$ such that
$\bar{||} \cdot \bar{||} \leq c \tilde{||} \cdot \tilde{||}$.  
At points where $G$ is not branched we can make a local 
expression $G = a z + b z^2 + \ldots$ with $a \neq 0$.  
We may assume $d\tilde{s}^2$ is the Euclidean metric 
locally, so $|dG|_{d\tilde{s}^2} = a$ at the 
chosen point .  At points where $G$ is branched we can make a local 
expression $\bar{G} = a z^{m} + b z^{m+1} + \ldots$ with 
$a \neq 0$ and $m \in \bfZ$, $m \geq 2$.  In this case 
$|dG|_{d\tilde{s}^2} = 0$ at the 
chosen point.  At each end we can make a local 
expression $G = z^{\mu} (a + b z + c z^2 + \ldots)$ with 
$a \neq 0$ and $\mu 
\geq 1$.  In this case $|dG|_{d\tilde{s}^2} = 0$ at the 
chosen point if $\mu > 1$, and $|dG|_{d\tilde{s}^2} = a$ 
at the chosen point if $\mu = 1$.  In any case $|dG|_{d\tilde{s}^2}$ 
is bounded, and the existence of $c$ follows.
\end{remark}

\begin{remark}
Since $d\bar{s}^2$ is identically zero for the horosphere, the 
calculations in this section would have no meaning for this example.  
And as it is the only example for which 
$d\bar{s}^2$ is zero at more than just isolated points, it is 
natural to exclude it.  In any case, the 
index of the horosphere is easily seen to be $0$ (see section 6).  
\end{remark}

\begin{lemma}
$\bar{H}^{1}$ is compactly contained in 
$L^{2}_{d\bar{s}^2}$.
\end{lemma}

\begin{proof}
$\bar{H}^{1}$ is continuously contained in 
$\tilde{H}^{1}$, 
and $\tilde{H}^{1}$ is compactly contained in 
$L^{p}_{d\tilde{s}^2}$ 
for any value of $p$, and $L^{p}_{d\tilde{s}^2}$ is continuously 
contained in $L^{2}_{d\bar{s}^2}$ if $p$ is large enough.  The 
composition of a continuous map 
and a compact map and a continuous map is compact.
\end{proof}

We remark that this section above and Lemma 5.2 below have an indirect, but 
close, relationship with the works of Troyanov and others on Riemannian surfaces 
with conical singularities \cite{HT}, \cite{T}.  

\section{The relationship between Ind($M$) and eigenvalues 
of $\bar{L}$}

The last lemma in the previous section will 
lead us to an argument that Ind($\Sigma$) is equal to 
the number of negative eigenvalues of $\bar{L}$ 
on $\Sigma$.  (We are about to define Ind($\Sigma$) and $\bar{L}$.)  
First we show that Rayleigh quotient ${\cal Q}$ (as defined in the next 
lemma) is well defined for any smooth function on 
$\Sigma$.  This next lemma 
will allow us to start the minimization process 
(i.e. ${\cal Q}(u) < \infty$ for some $u$) in the proof of the 
lemma which comes after it.

Before considering the next lemma, we define the relevant Jacobi 
operators.  
The original Jacobi operator on $M$ is 
$L = \triangle - |dG|_{ds^2}^{2} = \triangle + 
2 K$ on $\Sigma$.   The Jacobi operator created by pulling back 
the metric on the sphere via $G$ is $\bar{L} = 
\frac{-1}{K} \triangle - |dG|_{d\bar{s}^2}^{2} = 
\bar{\triangle} - 2$ on 
$\Sigma$.  Note that $\bar{L}$ is defined everywhere on $\Sigma$ 
except at the isolated 
points where $d\bar{s}^2 = 0$ and possibly at points that 
represent the ends of $M$.  The operator associated to the regular 
metric $d\tilde{s}^2$ is 
$\tilde{L} = \tilde{\triangle} - 
|dG|_{d\tilde{s}^2}^{2}$, and is defined on all of 
$\Sigma$, except possibly at points that 
represent the ends of $M$.  We have 
\[ \int_{\Sigma} u \tilde{L} u d\tilde{A} = 
\int_{\Sigma} u L u dA = 
\int_{\Sigma} u \bar{L} u d\bar{A} \; . \]
Since $L$, $\bar{L}$, and $\tilde{L}$ 
are not well defined only at isolated points of $\Sigma$, these 
integrals are well defined.  

\begin{lemma}
${\cal Q}(u) := \frac{\int_\Sigma u\bar{L}u d\bar{A}}
{\int_\Sigma u^2 d\bar{A}} 
< \infty$ for all $u \in C^{\infty}(\Sigma)$.
\end{lemma}

\begin{proof}
Locally at each end, we can normalize $G(z)$ to be $G = z^\mu 
(1 + a_1 z + a_2 z^2 + \ldots)$, where 
$z$ is contained in a neighborhood $U$ of 
$z=0$, and $z=0$ represents the end, and $d\tilde{s}^2$ is 
the Euclidean metric on $U$, and $\mu > 0$, and
$|dG|^{2}_{d\tilde{s}^2} \approx \frac{8 
\mu^{2}|z|^{2\mu-2}}{(1+|z|^{2\mu})^{2}}$.  (We showed this in the 
proof of Lemma 4.2.)  Since 
$\int_{\Sigma} u \tilde{L} u d\tilde{A} = 
\int_{\Sigma} u \bar{L} u d\bar{A}$, we can show that the 
numerator of ${\cal Q}(u)$ is finite 
by showing that $\int_{\Sigma} u \tilde{L} u d\tilde{A}$ is 
finite.  To show this, it is sufficient to show that 
$\int_{U} u \tilde{L} u d\tilde{A}$ is finite at each end, 
since $u \in C^{\infty}(\Sigma)$, and $d\tilde{s}^2$ 
is nonsingular on the compact $\Sigma$, and $\tilde{L}$ 
is nonsingular on $\Sigma$ away from the ends.  

Since $u \in C^{\infty}(\Sigma)$, we have that 
$u,u_{x},u_{y}$ are all bounded on $U$.  Since $d\tilde{s}^2$ is 
the Euclidean metric on $U$, we have that 
$d\tilde{A} = rdrd\theta$ in polar coordinates on $U$.  
Furthermore, we have $\mu>0$, hence 
\[ \left| \int_U \left( u_{x}^{2} + u_{y}^{2} - 
\frac{8 \mu^{2} r^{2\mu-2} u^2}{(1+r^{2\mu})^{2}} \right) 
\, r dr d\theta \right| < \infty \; . \]  
So the numerator of ${\cal Q}(u)$ is finite, and 
therefore $|{\cal Q}(u)| < \infty$.  
\end{proof}

Given any closed region $\Omega \subset \Sigma \setminus \{p_j\}$, we can 
consider the 
Dirichlet problem $Lu=\lambda u$ and $\bar{L}u=\lambda u$ on $\Omega$ 
with $u|_{\partial \Omega} = 0$.  In general, 
$L$ and $\bar{L}$ will have different eigenvalues on $\Omega$; 
however, supposing that $V$ is some vector space of functions with 
compact support on $\Omega$, ${\cal Q}(u) < 0$ for 
all $u \in V$ if and only if 
$\frac{\int_\Omega uLu dA}
{\int_\Omega u^2 dA} < 0$ for all $u \in V$.  We define 
Ind($\bar{L},\Omega$) to be the maximum possible dimension of a 
subspace of  functions in 
$C_0^\infty(\Omega)$ on which ${\cal Q}(u) < 0$.  We define 
Ind($L,\Omega$) to be the maximum possible dimension of a 
subspace of functions in 
$C_0^\infty(\Omega)$ on which $\frac{\int_\Omega uLu dA}
{\int_\Omega u^2 dA} < 0$.  
Thus Ind($\bar{L},\Omega$) = Ind($L,\Omega$).  We consider a 
sequence of regions $\{\Omega_{i}\}_{i=1}^\infty$ such that $\Omega_{i} 
\subset \Omega_{i+1}$ and $\cup \Omega_{i} = 
\Sigma \setminus \{p_j\}$.  We define 
Ind($\bar{L},M$) = $\lim_{i \to \infty}$ 
Ind($\bar{L},\Omega_i$), and we define 
Ind($L,M$) = $\lim_{i \to \infty}$ 
Ind($L,\Omega_i$).  It follows that 
Ind($\bar{L},M$) = Ind($L,M$).  And by the definition given in the 
second section, Ind$_u$($M$) = Ind($L,M$).  Defining 
Ind($\Sigma$) := 
Ind($\bar{L},\Sigma$) to be the maximum possible dimension of a 
subspace of functions in 
$C^\infty(\Sigma)$ on which ${\cal Q}(u) < 0$, we have 
\[ \mbox{Ind}_u(M) \leq \mbox{Ind}(\Sigma) \; .  \]

In order to explicitly compute Ind$_u$($M$) and to show that 
Ind$_u$($M$) = Ind($\Sigma$), we would like to know that 
Ind($\Sigma$) {\em equals the number of negative eigenvalues of} 
$\bar{L}$ {\em on} $\Sigma$.  That this holds (Corollary 5.1) 
can be concluded from the next lemma.  The fact that Ind($\Sigma$) 
equals the number of negative eigenvalues of 
$\bar{L}$ on $\Sigma$ is very useful for making explicit estimates of 
Ind($M$), as we shall see.  

\begin{lemma}
We can find weak solutions 
$u \in \bar{H}^1$ of $\bar{L} u = \lambda u$ on $\Sigma$ so 
that the following hold:
\begin{itemize}
\item The set of eigenvalues consists of an infinite sequence 
\[ \lambda_{1} < \lambda_{2} < \ldots \rightarrow \infty \; . \]
\item Each eigenvalue has finite multiplicity and the eigenspaces (of 
weak solutions) corresponding to distinct eigenvalues are 
$L_{d\bar{s}^2}^{2}$ orthogonal.  
\item The direct sum of the eigenspaces is dense in 
$L_{d\bar{s}^2}^{2}$ 
for the $L_{d\bar{s}^2}^{2}$ norm.  
\item Any eigenfunction $u$ of $\lambda_j$ is contained in 
$C^\infty(\Sigma \setminus \{p_j\})$ and satisfies 
$\bar{L} u = \lambda_j u$ in the classical sense on 
$\Sigma \setminus \{p_j\}$.  
\item Any eigenfunction $u$ of $\lambda_j$ is contained in 
$C^0(\Sigma)$.  
\end{itemize}
\end{lemma}

\begin{proof}
The Rayleigh-Ritz quotient as defined in Lemma 5.1 is 
\[ {\cal Q}(u) := \frac{\int_{\Sigma} 
|du|^{2}_{d\bar{s}^2} - 2 u^2 d\bar{A}}
{\int_{\Sigma} u^2 d\bar{A}} 
= \frac{\int_{\Sigma} 
|du|^{2}_{d\tilde{s}^2} - |dG|_{d\tilde{s}^2}^{2} 
u^2 d\tilde{A}}
{\int_{\Sigma} \frac{1}{2} |dG|_{d\tilde{s}^2}^{2} u^{2} 
d\tilde{A}} \; , \; \; \; \; \; u \in 
\bar{H}^{1} \; . \]  
The denominator is the $L_{d\bar{s}^2}^{2}$ norm.  

The proof of the first three items 
follows by standard variational 
arguments, such as in the arguments on pages 55-59 of 
\cite{Be}.  The only difference between the proof of the lemma 
above and 
the proof in \cite{Be} is that elliptic regularity is used there 
to show 
that the eigenfunctions corresponding to the eigenvalues are classical 
solutions of the eigenvalue problem (on all of $\Sigma$).  
In our case we only conclude 
that we have weak solutions to the eigenvalue problem.  However, 
we can simply ignore the arguments where elliptic 
regularity is used, and the remaining arguments in 
\cite{Be} are sufficient to prove the 
first three items in the above lemma, so we shall not repeat the 
arguments here.  

We remark that in order to apply these standard 
variational arguments, it is crucial that we know 
that $\bar{H}^1$ is {\em compactly} included in $L^2_{d\bar{s}^2}$.  
This is why we were focusing on proving Lemma 4.4 in the previous 
section.  

We now turn to proving the last two items in the lemma.  
Suppose $u$ is a weak solution of 
$\bar{L} u = \lambda_j u$, so 
$\tilde{\triangle} u = q u$ in the weak sense, where 
$q = (1 + \frac{\lambda_{j}}{2}) |dG|^{2}_{d\tilde{s}^2}$.  
Since $q \in 
C^{\infty} (\Sigma \setminus \{p_j\})$, it follows 
from elliptic regularity (\cite{GT}, Corollary 8.11) that 
$u \in C^{\infty} (\Sigma \setminus \{p_j\})$
and satisfies $\bar{L} u = \lambda_j u$ in the 
classical sense on $\Sigma \setminus \{p_j\}$.  

Consider a small neighborhood $U_j \subset \Sigma$ of the point $p_j$ 
representing an end.  
If $\mu_j \geq 1$ at $p_j$, then $q \in C^{0}(U_j)$, 
and therefore any eigenfunction $u$ is contained in $C^0(U_j)$ 
(\cite{GT}, Theorem 8.8 and Corollary 7.11).  If $\mu_j < 1$ at 
$p_j$, we will see in the next three paragraphs that $u$ is still 
contained in $C^0(U_j)$.  

To show $u \in C^{0}(\Sigma)$, we only need to show $u \in 
C^{0}(U_j)$, since $u$ is $C^\infty$ away from the $p_j$.  First we 
state theorem 17.1.1 from 
\cite{H}.  Consider a linear operator of order $m$, $P(x,D) = 
\sum_{|\alpha| \leq m} a_{\alpha}(x)D^{\alpha}$ in an open set $X 
\subset \bfR^{n}$.  We may assume $X=U_j$, since we may choose 
$d\tilde{s}^2$ to be the standard Euclidean metric on $U_j$.  
In our case $P(x,D) = \tilde{\triangle}$, $m=n=2$.
Suppose that $P_{m}(0,D) = \sum_{|\alpha| 
= m} a_{\alpha}(0)D^{\alpha}$ is elliptic.  This is certainly true 
for $\tilde{\triangle}$.  
Suppose also that $a_{\alpha}$ in continuous when $|\alpha| = m$, and 
that for some $r \in (1,\infty)$, 
$a_{\alpha} \in L_{loc,d\tilde{s}^2}^{\frac{n}{m-|\alpha|}}(U_j)$ 
if $m-|\alpha| < 
\frac{n}{r}$, and $a_{\alpha} \in 
L_{loc,d\tilde{s}^2}^{r+\epsilon}(U_j)$ for some 
$\epsilon > 0$ if $m-|\alpha| = \frac{n}{r}$, and $a_{\alpha} 
\in L_{loc,d\tilde{s}^2}^{r}(U_j)$ if $m-|\alpha| > \frac{n}{r}$.  
In our case all of the coefficients 
are constant, so these conditions will hold.  The theorem 
says that if all these conditions are satisfied and $U_j$ is a 
sufficiently small neighborhood of $p_j$, then there is a linear 
operator $E$ in $L_{d\tilde{s}^2}^{r}(U_j)$ such that 
\begin{itemize}
\item $f \in L_{d\tilde{s}^2}^{r}(U_j) \rightarrow D^{\alpha}Ef \in 
L_{d\tilde{s}^2}^{s}(U_j)$ is continuous if $r \leq s \leq \infty$ 
and $\frac{1}{s} \geq 
\frac{1}{r} - (\frac{m-|\alpha|}{n})$ with strict inequality if 
$s=\infty$, 
\item $P(x,D) Ef = f$, $f \in L_{d\tilde{s}^2}^{r}(U_j)$, 
\item $EP(x,D)$v = v if v$\in C_{0}^{\infty}(U_j)$.  
\end{itemize} 
In our case we will have 
$Ef = u$ and $f = q u$.  We will choose $1 < r = s \approx 1$.  

We now show that $q u \in L_{d\tilde{s}^2}^{r}(U_j)$ if $r$ is 
sufficiently close to 1.  Since $u \in \bar{H}^1$, Lemma 4.3 implies 
that $u \in \tilde{H}^1$.  And since $\tilde{H}^1$ is compactly 
contained in $L_{d\tilde{s}^2}^{p}$ for all $p \geq 2$, we have that 
$\int_{U_j} u^p d\tilde{A}$ is finite for all $p \geq 2$.  As we saw 
in the proof of Lemma 4.2, $\int_{U_j} q^t d\tilde{A}$ is finite if 
$t$ ($t>1$) 
is sufficiently close to 1.  Choose such a $t$ sufficiently close to 
1, and choose $p$ sufficiently large so that $\frac{1}{p} + 
\frac{1}{t} < 1$.  Define $r>1$, $r \approx 1$ by 
$\frac{1}{r} = \frac{1}{p} + \frac{1}{t}$.  By the Holder inequality, 
$(\int_{U_j} (q u)^{r} d\tilde{A} )^{\frac{1}{r}} \leq 
(\int_{U_j} q^{t} d\tilde{A} )^{\frac{1}{t}} 
(\int_{U_j} u^{p} d\tilde{A} )^{\frac{1}{p}} < \infty$.  

So we have that $q u \in L_{d\tilde{s}^2}^{r}(U_j)$, and, by the 
theorem from \cite{H} stated above, we have that 
the map $f \rightarrow D^{\alpha} Ef$ from 
$L_{d\tilde{s}^2}^{r}(U_j)$ to $L_{d\tilde{s}^2}^{s}(U_j)$ is 
continuous when $|\alpha| \leq 2$ and $s=r$.  Therefore, 
letting $Ef = u$ and $f = q u$, we have 
$||D^{\alpha}u||_{L_{d\tilde{s}^2}^{s=r} (U_j)} \leq c 
|| q u ||_{L_{d\tilde{s}^2}^r (U_j)} < \infty$ when $|\alpha| \leq 2$.  
This implies that $u$ is contained in the 
Sobolev space $W^{2,r}(U_j)$.  
Then, since $u \in W^{2,r}(U_j)$, we have $u \in C^{0}(U_j)$ 
(see, for example, \cite{GT}, Corollary 7.11).  
\end{proof}

\begin{corollary}
Ind($\Sigma$) = the number of negative eigenvalues (counted with 
multiplicity) of $\bar{L}$ on $\Sigma$.  
\end{corollary}

This corollary follows immediately from the variational 
characterization of the eigenvalues (for example, see \cite{Be}, 
p61).   The $k$'th eigenvalue $\lambda_k$ is characterized by 
\[ \lambda_k = \inf_{V_k} \sup \{Q(u) \, | \, u \in V_k, u 
\neq 0 \} \; , \] where $V_k$ runs through all $k$ dimensional 
subspaces of $\bar{H}^1$.  

\begin{lemma}
Ind$_u$($M$) = Ind($\Sigma$), and either 
Ind($M$) = Ind($\Sigma$) or Ind($M$) = Ind($\Sigma$) $-1$.
\end{lemma}

\begin{proof}
Let $f_1,\ldots,f_{\mbox{Ind}(\Sigma)}$ be the eigenfunctions 
of $\bar{L}$ on $\Sigma$ with negative eigenvalues.  Since $f_{i} \in 
\bar{H}^1$, we know that 
$\int_{B_\epsilon (p_j)} |f_{i}|^{2} d\bar{A} \rightarrow 0$ and 
$\int_{B_\epsilon (p_j)} |\bar{\nabla} f_{i}|^{2} d\bar{A} 
\rightarrow 0$ as $\epsilon \rightarrow 0$, where 
$B_\epsilon (p_j)$ is a 
ball about the end $p_j \in \Sigma$ with radius $\epsilon$ with 
respect to the $d\tilde{s}^2$ metric (these follow from the Lebesque 
dominated convergence theorem).  
By the previous lemma, we have $f_{i} \in C^{0}(\Sigma)$, and thus 
$|f_{i}| \leq c$, a constant.  Using these facts, we can then 
follow, with only slight modification, 
the argument in Fischer-Colbrie's proof \cite{FC}.  For the sake of 
completeness, we include the argument here.  

In a neighborhood of a point $p_j$ representing an end, choose a local 
complex coordinate $z$ centered at $p_j$.  For some small $\epsilon>0$, 
define a function 
$\eta_j(z) = 0$ if $|z| < \epsilon^2$, 
$\eta_j(z) = 1$ if $|z| > \epsilon$, and 
$\eta_j(z) = \frac{\log(\frac{|z|}{\epsilon^2})}
{\log(\frac{1}{\epsilon})}$ if $\epsilon^2 \leq |z| \leq \epsilon$.  
Let $\eta = \eta_j$ in an $\epsilon$ ball about each $p_j$, and let 
$\eta = 1$ elsewhere.  One can check that $\int_{\Sigma} 
|\bar{\nabla} \eta|^2 d\bar{A}  = \int_{\Sigma} 
|\tilde{\nabla} \eta|^2 d\tilde{A} \leq 
\frac{\hat{c}}{\log(\frac{1}{\epsilon})}$ for some constant $\hat{c}$, 
by noting that $d\tilde{s}^2 \approx |dz|^2$.  
Therefore $\int_{\Sigma} |\bar{\nabla} \eta|^2 d\bar{A} \to 0$
as $\epsilon \to 0$.  

Let $g_i = \eta f_i$, then 
\[ \int_{\Sigma} (g_i - f_i)^2 d\bar{A} = 
\int_{\Sigma} (1 - \eta)^2 f_i^2 d\bar{A} \leq 
\sum_{p_j} \int_{B_\epsilon(p_j)} |f_{i}|^{2} d\bar{A} \to 0 \]
as $\epsilon \to 0$, 
so $\int_\Sigma (g_i - f_i)^2 d\bar{A} \to 0$ as $\epsilon \to 
0$.  Also, 
\[ \int_\Sigma |\bar{\nabla} (g_i - f_i)|^2 d\bar{A} = 
\int_\Sigma |\bar{\nabla} ((1 - \eta)f_i)|^2 d\bar{A} \leq 
2 \int_\Sigma [|\bar{\nabla} \eta|^2 f_i^2 + 
(1 - \eta)^2 |\bar{\nabla} f_i|^2] d\bar{A} \leq \]\[
2c^2 \int_\Sigma |\bar{\nabla} \eta|^2 d\bar{A} + 
\sum_{p_j} \int_{B_\epsilon (p_j)} 
|\bar{\nabla} f_i|^2 d\bar{A} \; \; , 
\] since $|f_i| \leq c$ and $0 \leq (1-\eta)^2 \leq 1$.  Each of the 
integrals in the sum on the right hand side converge to 0 as 
$\epsilon \to 0$.  Hence we have that 
$\bar{||} f_i - g_i \bar{||}^2 \to 0$ as 
$\epsilon \to 0$.  By continuity of 
$\cal Q$ with respect to the $\bar{H}^1$ norm, we have that 
$\cal Q$ is negative definite 
on the span of $\{g_i\}_{i=1}^{\mbox{Ind}(\Sigma)}$ in $\bar{H}^1$ 
for sufficiently small $\epsilon$.  
Therefore Ind$_u$($M$) $\geq$ Ind($\Sigma$), and hence the first 
part of the lemma follows.

To prove the second part of the lemma, suppose that $V \subset 
\bar{H}^1$ is a vector space of dimension Ind$_u$($M$) on which 
${\cal Q}<0$.  If $V$ is perpendicular to the constant functions 
with respect to the $L^2_{d\bar{s}^2}$ inner product, 
then all of the functions in $V$ are volume preserving, and 
we have Ind($M$) = Ind($\Sigma$).  

If $V$ is not perpendicular to the constant functions, then the 
perpendicular projection of the constant function 1 to $V$ is a function 
$\phi_1 \neq 0, \phi_1 \in V$.  We may extend $\phi_1$ to an 
orthogonal basis $\{\phi_1, \phi_2,\ldots,\phi_{\mbox{Ind}_u(M)}\}$ 
of $V$.  Since $\phi_2,\ldots,\phi_{\mbox{Ind}_u(M)}$ are all 
perpendicular to $\phi_1$ with respect to the 
$L^2_{d\bar{s}^2}$ inner product in $V$, it 
follows easily that $\phi_2,\ldots,\phi_{\mbox{Ind}_u(M)}$ are all 
perpendicular to the constant function 1 
in $L^2_{d\bar{s}^2}$.  Thus a subspace of $V$ of dimension 
Ind$_u$($M$)$-1$ 
is perpendicular to the constant functions, so we have constructed a 
space of volume preserving functions of dimension Ind$_u$($M$)$-1$ on 
which ${\cal Q}<0$, and thus Ind($M$) $\geq$ Ind($\Sigma$)$-1$.  
\end{proof}

We have the following corollary, which is a result 
of do Carmo and Silveira (\cite{CS}).  The advantage of our proof of 
this corollary is that our 
method will allow us to make specific estimates of the index, whereas 
the method in \cite{CS} would not allow this.

\begin{corollary}
If a constant mean curvature 1 surface in $\bfH^3$ has 
finite total curvature, then it has 
finite index.
\end{corollary}

\begin{proof}
Since the surface has finite total curvature, it has a conformal 
bijection to $\Sigma \setminus \{p_j \}$.  By Lemma 5.2, $\bar{L}$ 
has a finite number of negative eigenvalues.  Then, by Corollary 5.1, 
Ind($\Sigma$) is finite.  Hence, by Lemma 5.3, Ind($M$) is finite.
\end{proof}

\section{Examples}

We now compute the index of several examples, showing how the results 
of the previous section can be applied.  For the sake of completeness, 
we compute the already known index of the horosphere, before 
continuing on to new results about index of certain surfaces.  

{\bf Horosphere:} For the horosphere, we can choose $\Sigma = 
\bfC \cup \{\infty\}$ and $\Sigma \setminus \{p_j\} = \bfC$ and 
$f=1$ and $g=1$ and $c=1$ and $z_0 = 0$ in Lemma 3.2.  Writing $F$ as 
\[ F = \left( \begin{array}{cc} 
A & B \\ 
C & D 
\end{array} \right) \; , \] we have 
\[ 
\left( \begin{array}{cc} 
A^\prime & B^\prime \\ 
C^\prime & D^\prime 
\end{array} \right) = 
\left( \begin{array}{cc} 
A & B \\ 
C & D 
\end{array} \right)
\left( \begin{array}{cc} 
1 & -1 \\ 
1 & -1 
\end{array} \right) \rightarrow
\left( \begin{array}{cc} 
A & B \\ 
C & D 
\end{array} \right) = 
\left( \begin{array}{cc} 
A_0 & B_0 \\ 
C_0 & D_0 
\end{array} \right) 
e^{
\left( \begin{array}{cc} 
1 & -1 \\ 
1 & -1 
\end{array} \right) z} \; \; \; , \]
and since $F |_{z_0 =0} = $id., we have $A_0 = D_0 = 1$ and $B_0 = C_0 = 
0$, therefore 
\[ 
F = \left( \begin{array}{cc} 
A & B \\ 
C & D 
\end{array} \right) = 
\left( \begin{array}{cc} 
1 & 0 \\ 
0 & 1 
\end{array} \right) + 
\left( \begin{array}{cc} 
1 & -1 \\ 
1 & -1 
\end{array} \right) z + 
\left\{ \left( \begin{array}{cc} 
1 & -1 \\ 
1 & -1 
\end{array} \right)^2 = 
\left( \begin{array}{cc} 
0 & 0 \\ 
0 & 0 
\end{array} \right) \right\}
\frac{z^2}{2!} + \ldots \]
and so 
\[ F = \left( \begin{array}{cc} 
1+z & -z \\ 
z & 1-z 
\end{array} \right) \; \; . \]  
Thus $G=1$, and it follows that the curvature 
$K=0$, and so the second variational formula given in the first 
section becomes 
\[ \left. \frac{d^2A}{dt^2} \right|_{t=0} = 
\int_M |\nabla u|^2 dA \geq 0 \; . \]  This is nonnegative for all 
functions $u$, hence the horosphere is stable.  

Silveira (\cite{Si}) showed that the only complete stable noncompact 
constant mean curvature 1
surface in $\bfH^3$ is the horosphere.  

{\bf Enneper cousin:} For the Enneper cousin, we can choose $\Sigma = 
\bfC \cup \{\infty\}$ and $\Sigma \setminus \{p_j\} = \bfC$ and 
$f=1$ and $g=z$ and $c=1$ and $z_0 = 0$ in Lemma 3.2.  Solving the 
equation 
\[ dF = F
\left( \begin{array}{cc} 
z & -z^2 \\ 
1 & -z  
\end{array} \right) dz \; \; , \]
we find that 
\[ F = 
\left( \begin{array}{cc} 
\cosh(z) & \sinh(z) - z \cosh(z) \\ 
\sinh(z) & \cosh(z) - z \sinh(z)  
\end{array} \right) \; . \] 
Therefore $G = 
\frac{d(\cosh(z))}{d(\sinh(z))} = \tanh(z)$.  

Following the Weierstrass representation as formulated in \cite{By}, 
we have a constant mean curvature 1 surface given by $F \bar{F}^t$ 
with {\em secondary} Gauss map $g=z$.  (Note that, since we are using 
$F$ instead of $F^{-1}$ to make the surface, the function $g$ is now the 
secondary Gauss map, not the hyperbolic Gauss map.)  In this case the 
secondary Gauss map is actually single valued, since the surface is 
simply connected.  By Lemma 4.1, the second variation is determined by 
$\Sigma = \bfC \cup \{\infty\}$ and $g=z$.  For this $\Sigma$ and 
$g$, the unconstrained index is Ind$_u$($M$) = 1.  This can be seen 
from Theorem 4.6 of \cite{N1}, or from 
Proposition 6.1 below.  It follows Ind($M$) is 
either 0 or 1.  But the Enneper cousin cannot be stable, since the 
horosphere is the only stable example (\cite{Si}), hence 
Ind($M$) = 1.  

We can also consider Enneper cousins with winding order $2k+1$ at the
end, $k \in \bfN$.  In this case $g$ becomes $g=z^k$, and the other 
objects $\Sigma = 
\bfC \cup \{\infty\}$ and $\Sigma \setminus \{p_j\} = \bfC$ and 
$f=1$ and $c=1$ and $z_0 = 0$ remain unchanged.  Now, by \cite{N1} or 
Proposition 6.1 below, Ind$_u$($M$) = $2k-1$.  Hence, by Lemma 5.3, the 
Enneper cousins with winding order $2k+1$ 
have constrained index Ind($M$) either $2k-1$ or $2k-2$.  

Following the Weierstrass representation as formulated in Lemma 3.2 
of this paper, 
we have a constant mean curvature 1 surface given by $F^{-1} 
\overline{F^{-1}}^t$, and this produces the ``dual'' Enneper cousin 
that is described in \cite{RUY}.  The dual Enneper cousin has 
secondary Gauss map $G = \tanh(z)$ 
and hence has infinite total curvature.  By \cite{CS}, it must 
therefore have infinite index.  (See Figure 1.)  

\begin{figure}
        \hspace{2in}
        \epsfxsize=1.8in
        \epsffile{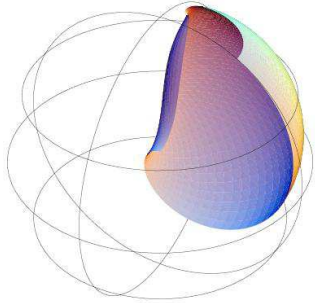}
        \hfill
\caption{Half of a ``dual'' Enneper cousin in the Poincare model.
The entire surface consists of the piece shown union its reflections 
across planes containing boundary curves.  This 
surface has infinite total curvature, and therefore has infinite 
index.}
\end{figure}

{\bf Catenoid cousins:} A catenoid cousin has $\Sigma \setminus 
\{p_j\} = \bfC \setminus \{0\}$, and has secondary Gauss map $G = z^\mu$,
where $\mu \neq 0, \pm 1$ is real.  We can assume without loss of
generality that $\mu > 0$.  The surface is embedded if $\mu <1$ and
not embedded if $\mu > 1$.  (This is shown in \cite{UY1}.  It 
was originally shown in \cite{By}, but the parameter $\mu$ is 
formulated differently in \cite{By}.  We use the same $\mu$ as in 
the \cite{UY1} formulation.  Figures of the catenoid cousins can 
be found in \cite{UY1}.)  

We will show that the embedded catenoid cousins have index 1, and that 
the non-embedded catenoid cousins have index at least 2, and that the 
index gets arbitrarily large as $\mu$ gets arbitrarily large.  To 
prove this, we first prove the following proposition.  
This proposition is proven in \cite{N1} in the case that
$\mu$ is an integer.  The proof when $\mu$ is not an integer is 
essentially the same.  We include the proof here for the sake of 
completeness.

\begin{proposition}
Let $\mu$ be a positive real.  
The complete set of eigenvalues for the Laplacian on the plane with 
the pull back metric from the sphere via the map $G = z^\mu$ is 
\[ \lambda_{p,q} = (p+\frac{q}{\mu})(1 + p+\frac{q}{\mu}) \; , \; 
\; p,q=0,1,2,... \; \; \; \; \; . \]
The multiplicity of $\lambda_{p,q}$ is 2 if $q>0$ and 
is 1 if $q=0$.
\end{proposition}

\begin{proof}
We are considering the problem $\bar{\triangle} u = \lambda u$, 
where $\bar{\triangle}$ is the Laplacian obtained from pulling 
back the standard metric on $S^2$ via the map $G = z^\mu$.  In polar 
coordinates this equation becomes
\[ \frac{\partial^2 u}{\partial r^2} + 
\frac{1}{r} \frac{\partial u}{\partial r} + 
\frac{1}{r^2} \frac{\partial^2 u}{\partial \theta^2} = - 
\lambda \frac{4 \mu^2 r^{2\mu-2}}{(r^{2\mu} + 1)^2} u \; \; . \]

For a real number $\alpha$ and a nonnegative integer $i$, we define 
$(\alpha)_i$ to be $(\alpha)_i = \alpha(\alpha+1).....(\alpha+i-1)$ 
if 
$i > 0$ and $(\alpha)_i = 1$ if $i=0$.  We then define a real 
analytic 
hypergeometric function $F(a,b,c,x)$ where $c$ is not a nonpositive 
integer.  This function $F(a,b,c,x)$ is defined for $-1<x<1$.
\[ F(a,b,c,x) := \sum_{i=0}^\infty \frac{(a)_i (b)_i}{i! (c)_i} x^i 
\; \; \; . \]
$F(a,b,c,x)$ satisfies the hypergeometric differential equation 
\[ x(1-x)\frac{d^2y}{dx^2}+(c-(a+b+1)x)\frac{dy}{dx} - aby = 0 
\; \; \; . \]
For nonnegative integers $p$ and $q$, $F(p + 2 \frac{q}{\mu} + 1,
-p, \frac{q}{\mu} + 1, \frac{1}{2}(1-t))$ is a polynomial of degree 
$p$.
We set 
\[ \phi_{p,q}(t) = (1-t^2)^{\frac{q}{2\mu}} F\left(p + 
2 \frac{q}{\mu} + 1,
-p, \frac{q}{\mu} + 1, \frac{1}{2}(1-t)\right), \; \; -1<t<1 \]
and 
\[ v_{p,q}(r) = 
\phi_{p,q}\left(\frac{r^{2\mu}-1}{r^{2\mu}+1}\right), \; \; \; 
0<r<\infty \; \; \; . \]
We can check that, for $\lambda = (p+\frac{q}{\mu})(1 + 
p+\frac{q}{\mu})$, $\phi_{p,q}(t)$ satisfies the ordinary differential
equation
\[ (1-t^2)\frac{\partial^2 \phi}{\partial t^2}-2t\frac{\partial \phi}
{\partial t} + \left(\lambda - (\frac{q}{\mu})^2 
\frac{1}{1-t^2}\right)\phi = 0 
\; \; \; , \]
and $v_{p,q}(r)$ satisfies the ordinary differential equation
\[ \frac{\partial^2 v}{\partial r^2} + \frac{1}{r} 
\frac{\partial v}{\partial r} + \left(\lambda 
\frac{4 \mu^2 r^{2\mu-2}}{(r^{2\mu}+1)^2} - 
\frac{q^2}{r^2}\right) v = 0 \; \; \; . \]
We can then check that $v_{p,q}(r) \cos(q\theta)$ and 
$v_{p,q}(r) \sin(q\theta)$ are eigenfunctions of the Laplacian 
with eigenvalue $\lambda = (p+\frac{q}{\mu})(1 + p+\frac{q}{\mu})$.

Finally, we need to check that, for nonnegative $p$ and $q$,
the above eigenfunctions form a complete orthogonal system in the
$L^2_{d\bar{s}^2}$ norm.  This follows by elementary arguments.  
%
\end{proof}

The following theorem follows immediately from 
Proposition 6.1 and Lemma 5.3, and from 
Silveira's result that the horosphere is the only 
stable complete constant mean curvature 1 surface in 
$\bfH^3$ \cite{Si}.  

\begin{theorem}
The index of any embedded catenoid cousin is exactly 1, and the 
index of any nonembedded catenoid cousin is at least 2.  Let $[\mu]$ 
be the greatest integer that is strictly less than $\mu$, then the 
index of the catenoid with value $\mu$ is either $2[\mu]+1$ or 
$2[\mu]$.  Thus, 
for any positive number $N$ there exists a catenoid cousin with index
greater than $N$.
\end{theorem}

\begin{remark} It is clear from the above proposition that 
when $\mu \not\in \bfZ$, the nullity (nullity := the dimension of the 
eigenspace corresponding to the eigenvalue 0) of the catenoid cousins is 1 
(i.e. this is the case that p=1,q=0), and that when $\mu \in \bfZ$, the 
nullity is 3 ($p=1,q=0$ or $p=0,q=\mu$).  And the unconstrained 
index Ind$_u$($M$) changes only as $\mu$ passes through an integer, 
when two eigenvalues pass through 0.   
Furthermore, this illustrates another difference from the case of minimal 
surfaces in $\bfR^3$, where the nullity is always at least 3, since 
the set of translations make bounded normal Jacobi fields on a 
minimal surface (\cite{N1}, \cite{MR}, \cite{EK}).  
\end{remark}

There are some other examples where we can compute the 
index explicitly, via the above proposition, which we will now 
describe.  

Example 7.4 of \cite{UY1} has a Weierstrass 
representation with $G = z^{m}$, $3 \leq m \in 
\bfZ$ on $\Sigma \setminus \{p_j\} = \bfC \setminus \{0\}$.  It 
follows immediately from Proposition 6.1 and Lemma 5.3 that the index 
Ind($M$) of this example is either $2m-1$ or $2m-2$.  

Another example is given in 
Theorem 6.2 of \cite{UY1}. It has a Weierstrass representation with 
$G = a z^{\ell} + b$ and Hopf differential 
$Q = a c \ell z^{-2} (dz)^2$ on $\Sigma \setminus \{p_j\} = \bfC 
\setminus \{0\}$, where $\ell \in \bfZ, \ell \neq 0$ and 
$a,b,c \in \bfC$, $a \neq 0$, $c \neq 0$, and $\ell^{2} + 
4 a c \ell = m^{2}$ for some positive integer $m$.  

In the case that $b=0$, we can simply rewrite $a^{\frac{1}{\ell}} z$ 
as $z$, and then $\Sigma \setminus \{p_j\}$ is unchanged and $G$ 
becomes $G=z^\ell$.  Hence, when $b=0$, Ind($M$) is either $2\ell-1$ 
or $2\ell-2$, by Proposition 6.1 and Lemma 5.3.  

In the case that $\ell = 1$, then we can make the transformation 
of the complex plane $z \to \frac{z-b}{a}$.  Then $\Sigma$ is still 
$\bfC \cup \{\infty\}$, and $G$ becomes $G=z^1$.  Hence 
by Proposition 6.1 and Lemma 5.3, Ind($M$) is either 0 or 1.  
By \cite{Si} these surfaces cannot be stable, hence Ind($M$) = 1.
There are many different examples of this type 
with $\ell=1$: for example, $\ell=1,a=1,c=2,m=3,b=0$ or 
$\ell=1,a=\frac{3}{4},c=1,m=2,b=0$, and infinitely many others.

\begin{remark}
This last example with $\ell=1$ illustrates another difference between 
minimal surfaces in $\bfR^3$ and constant mean curvature surfaces 
in $\bfH^3$: 
While the only complete minimal surfaces in $\bfR^3$ with index 1 
are the catenoid and Enneper's surface (\cite{FC}, \cite{Cho}), the 
embedded catenoid cousins and the 
Enneper cousins are not the only constant mean curvature 1 surfaces 
in $\bfH^3$ with index 1.  
\end{remark}

\section{Lower bounds for Ind($M$)}

Choe \cite{Cho} proved some general results about lower bounds for the 
index of minimal surfaces in $\bfR^3$.  In this section, we will apply 
the same method to constant mean curvature 1 surfaces in $\bfH^3$.  
The results in section 5, particularly Lemma 5.3, are crucial to 
getting the method to work in our situation.  

Let $\phi$ be a Killing vector field in $\bfH^{3}$ generated by 
either a hyperbolic rotation or a hyperbolic translation.  
For both a hyperbolic rotation and a hyperbolic 
translation there are two fixed points on the sphere at infinity, and 
we shall call these two points the points in the sphere at infinity 
fixed by $\phi$.  (For example, a Euclidean rotation about 
the $x_3$-axis of the upperhalf space model for $\bfH^3$ is a 
hyperbolic rotation, and a Euclidean dilation centered at the point 
$x_1 = x_2 = x_3 = 0$ of the upperhalf space model is a 
hyperbolic translation.  Both of these isometries of $\bfH^3$ fix the 
two points $x_1 = x_2 = x_3 = 0$ and $x_1 = x_2 = x_3 = \infty$ in the 
sphere at infinity.)  

Let $M$ be a constant mean curvature 1 surface in $\bfH^{3}$ with 
finite total curvature.  The 
Killing vector field $\phi$ can be decomposed into tangent and normal 
parts on $M$, that is, $\phi = \phi^T + \phi^\perp$, where $\phi^T \in 
T(M)$ and $\phi^\perp \in N(M)$, and $N(M)$ is the normal vector 
bundle of $M$.  Choosing a unit normal $\vec{N}$ on $M$, it 
can be checked by a direct computation (see, for example, 
Lemma 1 of \cite{Cho} or Proposition 2.12 of \cite{BCE}) that the 
normal projection $\phi^\perp = u \vec{N}$ of a 
Killing vector field $\phi$ on $M$ is a Jacobi 
field (i.e. $\triangle u + 2 K u = 0$).  

\begin{definition}
Let $H(M,\phi)$ be the set of all points on $M$ where $\phi^\perp = 
0$.  We call $H(M,\phi)$ the {\em horizon} of 
$M$ with respect to $\phi$.  Each component of $M \setminus 
H(M,\phi)$ is called a {\em visible set} of $M$.  The number of 
visible sets of $M \setminus 
H(M,\phi)$ is called the {\em vision number} $v(M,\phi)$ of $M$ with 
respect to $\phi$.  The number of visible sets of $M \setminus 
H(M,\phi)$ which are either bounded or whose closure intersects the 
sphere at infinity only at one or both of the points fixed by $\phi$ 
is called the {\em adjusted vision number} $\tilde{v}(M,\phi)$ of $M$ with 
respect to $\phi$.  
\end{definition}

Note that $\tilde{v}(M,\phi) \leq v(M,\phi)$.

\begin{theorem}
Let $M$ be a constant mean curvature 1 surface in $\bfH^{3}$ 
of finite total curvature 
with regular ends.  Then for any choice of $\phi$, 
\[ \mbox{Ind}(M) \geq \tilde{v}(M,\phi) - 1 \] if 
$\tilde{v}(M,\phi) \neq v(M,\phi)$, and 
\[ \mbox{Ind}(M) \geq \tilde{v}(M,\phi) - 2 \] if 
$\tilde{v}(M,\phi) = v(M,\phi)$.  
\end{theorem}

\begin{proof}
First we show that 
on any visible set which is counted in $\tilde{v}(M,\phi)$, $u \vec{N}$ 
is bounded.  (This is not true 
for the visible sets which are not counted in $\tilde{v}(M,\phi)$.)  
For this, we need to use that we have regular 
finite total curvature ends.  Note that the definition of a regular 
end is an end for which the hyperbolic Gauss map $G$ extends 
holomorphically across $p_j$ \cite{UY1}.  
An end with finite total curvature 
is regular if and only if ord$_{p_j}(Q) \geq -2$ \cite{By}.  
For these types of ends, assuming that the end 
approaches the origin in the upper-half-space model, we have the 
following asymptotic behavior:
\[ \left( \Re(z^{m}) , \Im(z^{m}) , c |z|^{\mu+m} 
(1 + {\cal O}(|z|^{\mbox{min}(1,2\mu)})) \right) \; , \] 
where $z$ is a local coordinate at the end, $z=0$ is the point 
representing the end, and $c$ is a positive constant.  
${\cal O}(1,2\mu) = {\cal O}(|z|^{\mbox{min}(1,2\mu)})$ denotes 
any real valued function $f(z)$ such that $\limsup_{z \to 0} 
\frac{f}{|z|^{\mbox{min}(1,2\mu)}}$ is bounded.  
(See the appendix for a proof of this asymptotic 
behavior.)  Note that the unit normal $\vec{N}$ is of the form 
\[ \vec{N} = \frac{cx^{\mu+m}}
{\sqrt{m^2+c^2 (\mu+m)^2 x^{2\mu} + {\cal O}(1,2\mu)}} (-c(\mu+m)x^\mu 
(1+{\cal O}(1,2\mu)),x^\mu {\cal O}(1,2\mu),m) \]
over a point $z = x > 0, x \in \bfR$.  

We now consider three cases for the Killing 
vector field $\phi$:
\begin{itemize}
\item Suppose $\phi$ is made by an isometry which is either 
a hyperbolic rotation or a hyperbolic translation, and suppose that the 
origin $x_1 = x_2 = x_3 = 0$ 
is not one of the two points in the sphere at infinity fixed 
by $\phi$.  In this case we may consider that 
$\phi \approx \vec{(1,0,0)}$ near the origin.  Thus, when $z=x>0$, we 
have \[ \langle \phi, \vec{N} \rangle_{\bfH^3} \approx 
\frac{-(\mu+m)(1+{\cal O}(1,2\mu))}{x^m \sqrt{m^2+
c^2 (\mu+m)^2 x^{2\mu} + {\cal O}(1,2\mu)}} \; , \] 
and this will diverge to $\infty$ as $x \rightarrow 0$.  Thus 
for a $\phi$ of this type, the normal Jacobi vector field 
$\phi^\perp = \langle \phi, \vec{N} \rangle \vec{N}$ is not bounded.  
(As a simple example, one can easily compute $\langle \phi, \vec{N} 
\rangle$ explicitly for a horosphere.)  
\item Suppose $\phi$ is made by the isometry which 
is a dilation centered at the origin.  In this case 
$\phi = \vec{(x_1,x_2,x_3)}$ at $(x_1,x_2,x_3)$.  Thus, when $z=x>0$, we 
have 
\[ \langle \phi, \vec{N} \rangle_{\bfH^3} = 
\frac{-\mu (1+{\cal O}(1,2\mu))}{\sqrt{m^2+c^2 
(\mu+m)^2 x^{2\mu} + {\cal O}(1,2\mu)}} \approx \frac{-\mu}{m}
\; . \]  Thus, for a $\phi$ of this type, the length $\langle \phi, 
\vec{N} \rangle_{\bfH^3}$ of the normal Jacobi 
vector field $\phi^\perp$ is bounded and continuous in a neighborhood 
of the end.  
\item Suppose $\phi$ is made by the isometry which is rotation about 
the $x_3$-axis.  In this case 
$\phi = \vec{(-x_2,x_1,0)}$ at $(x_1,x_2,x_3)$.  Thus, when $z=x>0$, we 
have \[ \langle \phi, \vec{N} \rangle_{\bfH^3} = 
{\cal O}(1,2\mu) \; , \] and so the length 
$\langle \phi, \vec{N} \rangle_{\bfH^3}$ of the normal Jacobi 
vector field $\phi^\perp$ is bounded and continuous in a 
neighborhood of the end.  
\end{itemize}

Let $u = \langle \phi, \vec{N} \rangle_{\bfH^3}$ be the length 
of the normal variation vector field $\phi^\perp$.  In the second 
and third cases above, $u$ is bounded and continuous 
at the end asymptotic to the origin in the upper half space model.  
Hence we can conclude from Harvey and Polking's removable singularity 
theory (\cite{HP}, \cite{P}, \cite{Cho}) that $u$ is a weak solution of 
the Jacobi operator 
$\triangle u + 2 K u = 0$ on $\Sigma$, except at the ends where $u$ 
is not bounded.  
%

So on each visible set counted in $\tilde{v}(M,\phi)$, $u$ is 
bounded; and for each visible set counted in $\tilde{v}(M,\phi)$, 
the nullity of the visible set with respect to the Dirichlet 
problem is at least 1.  

The operator $\bar{L}$ on $\Sigma$ has the following properties:
\begin{itemize}
\item $\bar{L}$ satisfies the {\em unique continuation 
property}; that is, if two solutions $u$ and $v$ of $\bar{L} = 0$ 
are equal on any open 
set of $\Sigma$, then they are equal on all of $\Sigma$.  
This property holds on $\Sigma$ simply because it holds 
on $\Sigma \setminus \{ p_j \}$ (since any weak 
solution $u$ of $\bar{L} u = 0$ on $\Sigma \setminus \{ p_j \}$
is also a strong solution on $\Sigma \setminus \{ p_j \}$, by 
elliptic regularity), and because any open set in 
$\Sigma$ contains an open set of $\Sigma \setminus \{ p_j \}$.  
\item $\bar{L}$ also satisfies a 
{\em variational characterization of the eigenvalues property}.  This 
follows from the standard variational arguments used in the proof of 
the first three items of Lemma 5.2.  The eigenvalues can be 
characterized as $\lambda_k = \min({\cal Q}(u))$, where the functions 
$u \neq 0$ are any functions that are $L^2_{d\bar{s}^2}$-perpendicular 
to the eigenspaces of $\lambda_1,\ldots,\lambda_{k-1}$.  And those 
functions $u$ for which $\cal Q$ attains the minimum $\lambda_k$ are 
precisely the eigenfunctions associated to $\lambda_k$.  
\item Using the above two properties and the variational 
characterization we used to derive Corollary 5.1, we can conclude that 
as a domain $\Omega$ increases in size, the eigenvalues (with respect 
to the Dirichlet problem) must be strictly decreasing.  
\end{itemize}
These properties enable us to conclude that Smale's theorem holds in 
our setting \cite{FT}.  

{\em Smale's result}:  (see \cite{L}, Theorem 33) 
Let $c_t$ be a smooth contraction of $\Sigma$ into itself such that 
\begin{itemize}
\item $c_0$ = identity
\item $c_t(\Sigma) \subset c_s(\Sigma)$ for $t>s$
\item $\lim_{t \rightarrow \infty} \mbox{Volume}(c_t) = 0$
\end{itemize}
then 
\[ \mbox{Ind}(\Sigma) \geq \sum_{t>0} \mbox{Nullity} (c_t) \; , \] 
where $\mbox{Nullity} (c_t)$ is the dimension of the space of Jacobi 
fields on $c_t(\Sigma)$ vanishing on the boundary of 
$c_t(\Sigma)$.  

Noting that we have shown that Ind($\Sigma$) = Ind$_u$($M$), 
the proof then follows essentially as in the proofs of Theorem 1 
of \cite{Cho}.  For the sake of completeness, we include the argument 
here.  

Let $k = v(M,\phi)$ and $\tilde{k} = \tilde{v}(M,\phi)$, and let 
$V_1,\ldots,V_k$ be the open components of $M \setminus H(M,\phi)$.  
Let $\hat{V}_1,\ldots,\hat{V}_k$ be the open sets of $\Sigma$ 
corresponding to the sets $V_1,\ldots,V_k$ under the 
conformal bijection between $\Sigma$ and $M$.  After suitably 
renumbering $\hat{V}_1,\ldots,\hat{V}_k$, we can exhaust the sets 
$\hat{V}_j$ by a continuous 1-parameter 
family of shrinking domains $c_t(\Sigma),t \in (0,\infty)$ with 
piecewise smooth boundaries such that 
$\Sigma \setminus c_t(\Sigma) \subset \hat{V}_1$ for $t<1$, and 
$c_j(\Sigma) = 
\hat{V}_{j+1} \cup \ldots \cup \hat{V}_k$ for each integer $j = 
1,\ldots,k-1$, and $c_t(\Sigma) \subset \hat{V}_k$ for all $t>k-1$.  
We may assume that $V_{k-\tilde{k}+1},\ldots,V_k$ are the 
sets that are counted in $\tilde{v}(M,\phi)$, and that 
$V_1,\ldots,V_{k-\tilde{k}}$ are not counted in $\tilde{v}(M,\phi)$.

If $\tilde{v}(M,\phi) < v(M,\phi)$, it follows that 
Nullity$(c_t(\Sigma)) \geq 1$ with respect to the 
Dirichlet problem on $c_t(\Sigma)$ when $t=k-\tilde{k}, 
t=k-\tilde{k} + 1, \ldots, t = k-1$.  By Smale's theorem it follows 
that \[ \mbox{Ind}(\Sigma) \geq \tilde{k} \; \; \; . \]
If $\tilde{v}(M,\phi) = v(M,\phi)$, it follows that 
Nullity$(c_t(\Sigma)) \geq 1$ with respect to the 
Dirichlet problem on $c_t(\Sigma)$ when $t=1, 
t=2, \ldots, t = k-1$.  By Smale's theorem it follows 
that \[ \mbox{Ind}(\Sigma) \geq \tilde{k}-1 \; \; \; . \]
Thus, by Lemma 5.3, the theorem is proved.  
\end{proof}

We now apply Theorem 7.1 to find lower bounds for the index of 
several specific examples.  

\begin{corollary}
Suppose that $M$ is a constant mean curvature 1 genus $k$ Costa cousin 
in $\bfH^3$ (as described in \cite{RUY}, $M$ is the surface in $\bfH^3$ 
corresponding to the minimal genus $k$ Costa-Hoffman-Meeks surface).  
Then Ind($M$)$\geq 2k$.  
\end{corollary}

\begin{proof}
Consider the surface $M$ in the Poincare model with two ends 
asymptotic to the point $(0,0,1)$ in the sphere at infinity and one 
end asymptotic to the point $(0,0,-1)$ in the sphere at infinity (see 
figure 2).  Let $\phi$ be the Killing vector field generated by 
hyperbolic rotation about the $x_3$-axis, thus $\phi$ fixes the two 
points $(0,0,1)$ and $(0,0,-1)$ in the sphere at infinity.  
Due to the reflective symmetries of $M$, it is clear that 
$\tilde{v}(M,\phi) = v(M,\phi) \geq 2k+2$.  By Theorem 7.1, the 
corollary follows.  
\end{proof}

Using the same $\phi$ as in the above proof and placing the genus 1 
catenoid cousins (as described in \cite{RS}) so that their ends are 
asymptotic to $(0,0,1)$ and $(0,0,-1)$ (see figure 3), we also have 
the following corollary.

\begin{corollary}
The genus one catenoid cousins have index at least 2.  
\end{corollary}

The following proposition can be proven by an argument similar to that 
of Corollary 4 in \cite{Cho}.  This result is not stated in 
\cite{Cho}, perhaps 
only because existence of the minimal genus 1 $n$-noid was not 
known at that time \cite{BR}.

\begin{proposition}
The minimal genus 1 $n$-noid in $\bfR^3$ has index at least $n$ if 
$n$ is odd, and at least $n-1$ if $n$ is even.  
\end{proposition}

The corollaries above do not require that the surfaces be slight 
deformations of a corresponding minimal surface, but in the following 
corollary we will need this assumption.  The deformation is 
described in \cite{RUY}.  

\begin{corollary}
Let $M$ be a constant mean curvature 1 genus 1 $n$-noid cousin in 
$\bfH^3$ that is a slight deformation of a minimal genus 1 $n$-noid 
in $\bfR^3$ (see figure 4).  Then 
Ind($M$)$\geq n-3$ if $n$ is even and 
Ind($M$)$\geq n-4$ if $n$ is odd.  
\end{corollary}

\begin{proof}
Place the surface $M$ in the Poincare model so that all of its $n$ 
ends are asymptotic to points in the sphere at infinity where $x_2 = 
0$.  Let $p_1$ and $p_2$ be two points in the sphere at infinity such 
that two adjacent ends of $M$ are asymptotic to $p_1$ and $p_2$.  
Let $\phi$ be a hyperbolic translation fixing $p_1$ and $p_2$ in the 
sphere at infinity.  When $M$ is a sufficiently small deformation of a 
minimal genus 1 
$n$-noid, we know the behavior of $H(M,\phi)$ (since we know 
what the horizon is on the minimal genus 1 
$n$-noid \cite{Cho}).  (See figure 
5.)  We can conclude that $n-2 = \tilde{v}(M,\phi) < v(M,\phi)$ when 
$n$ is even, and $n-3 = \tilde{v}(M,\phi) < v(M,\phi)$ when $n$ is 
odd.  Then Theorem 7.1 implies the corollary.  
\end{proof}

\begin{figure}
        \hspace{1.1in}
        \epsfxsize=1.8in
        \epsffile{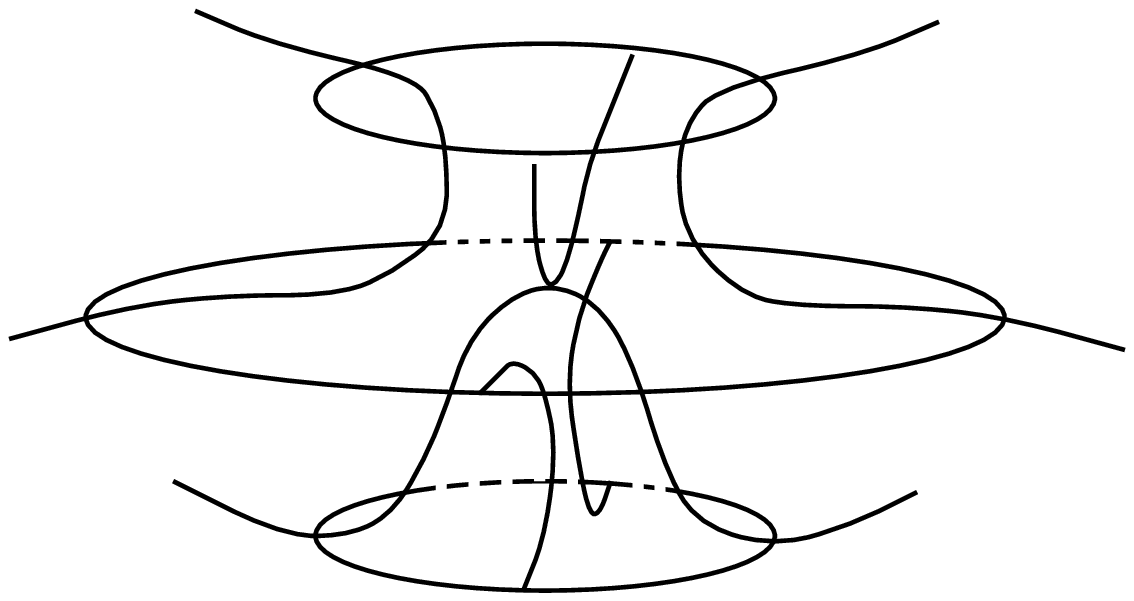}
        \hspace{0.01in}
        \epsfxsize=1.1in
        \epsffile{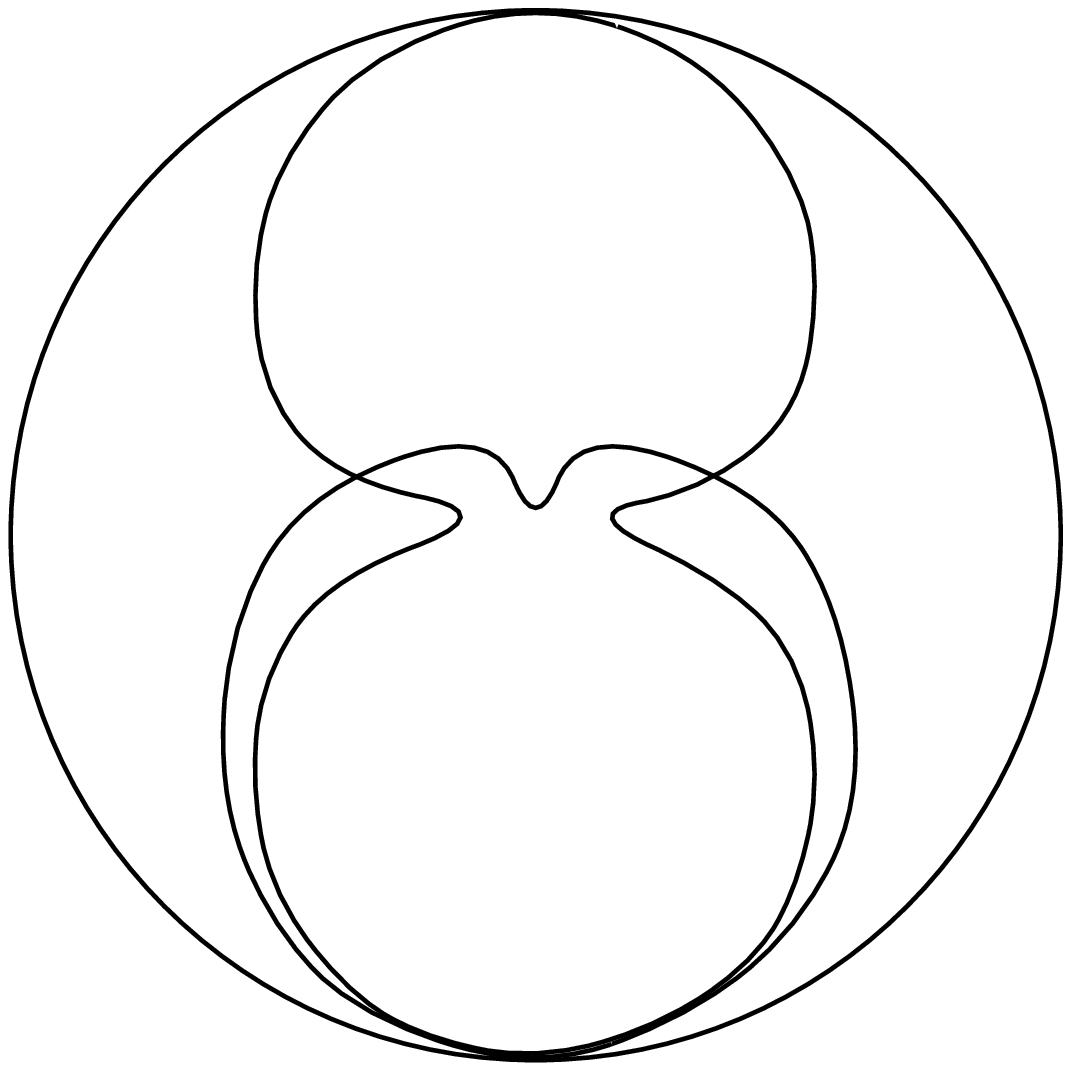}
        \hspace{0.01in}
        \epsfxsize=1.1in
        \epsffile{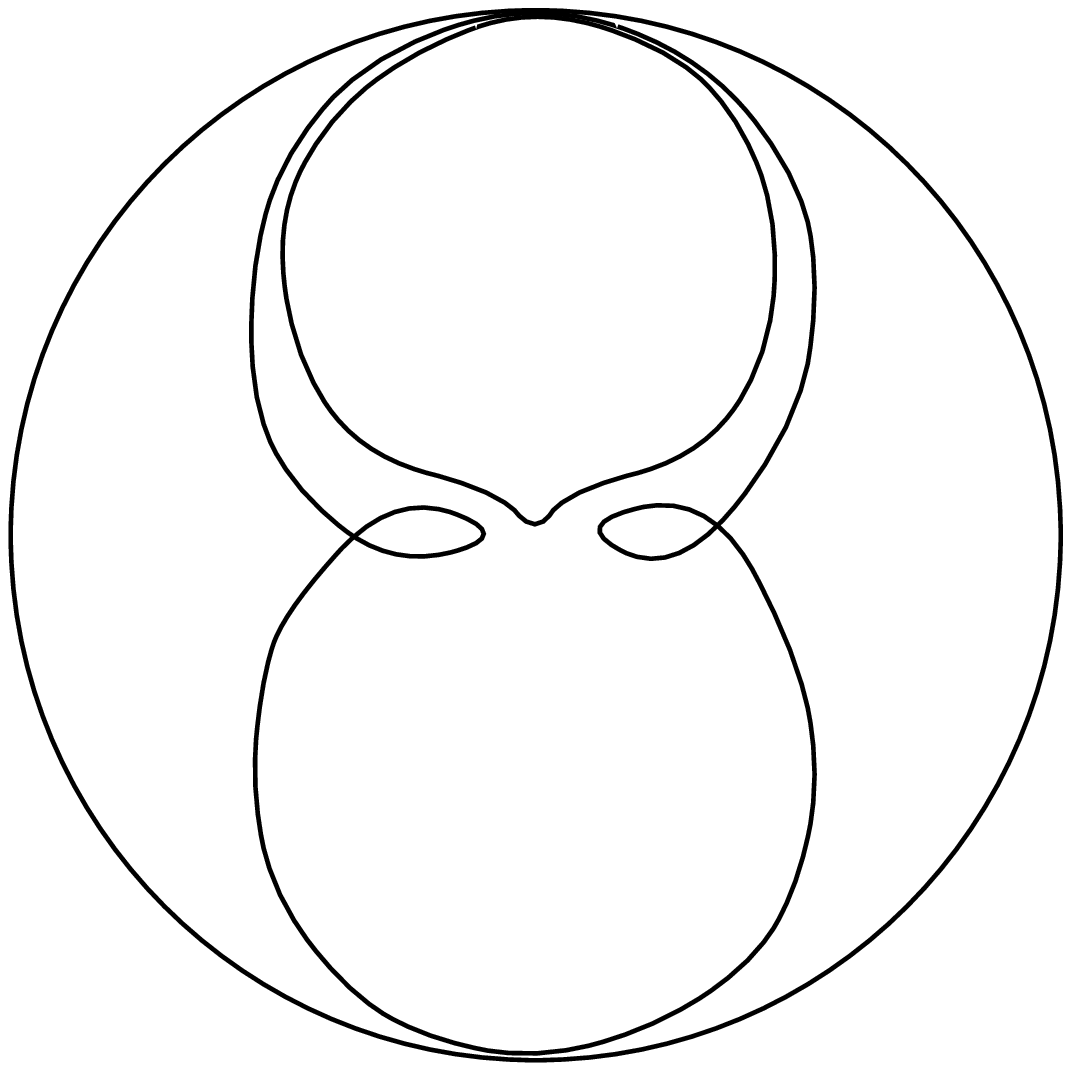}
        \hfill
\caption{A genus 1 Costa cousin in the Poincare model: slices in 
the $x_1 x_3$-plane and in the $x_2 x_3$-plane.}
\end{figure}

\begin{figure}
        \hspace{2in}
        \epsfxsize=1.8in
        \epsffile{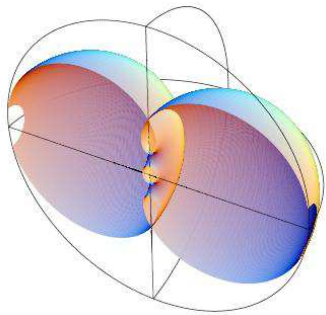}
        \hfill
\caption{Half of a genus 1 catenoid cousin (computer graphics by 
Katsunori Sato of Tokyo Institute of Technology).}
\end{figure}

\begin{figure}
        \hspace{1in}
        \epsfxsize=1.8in
        \epsffile{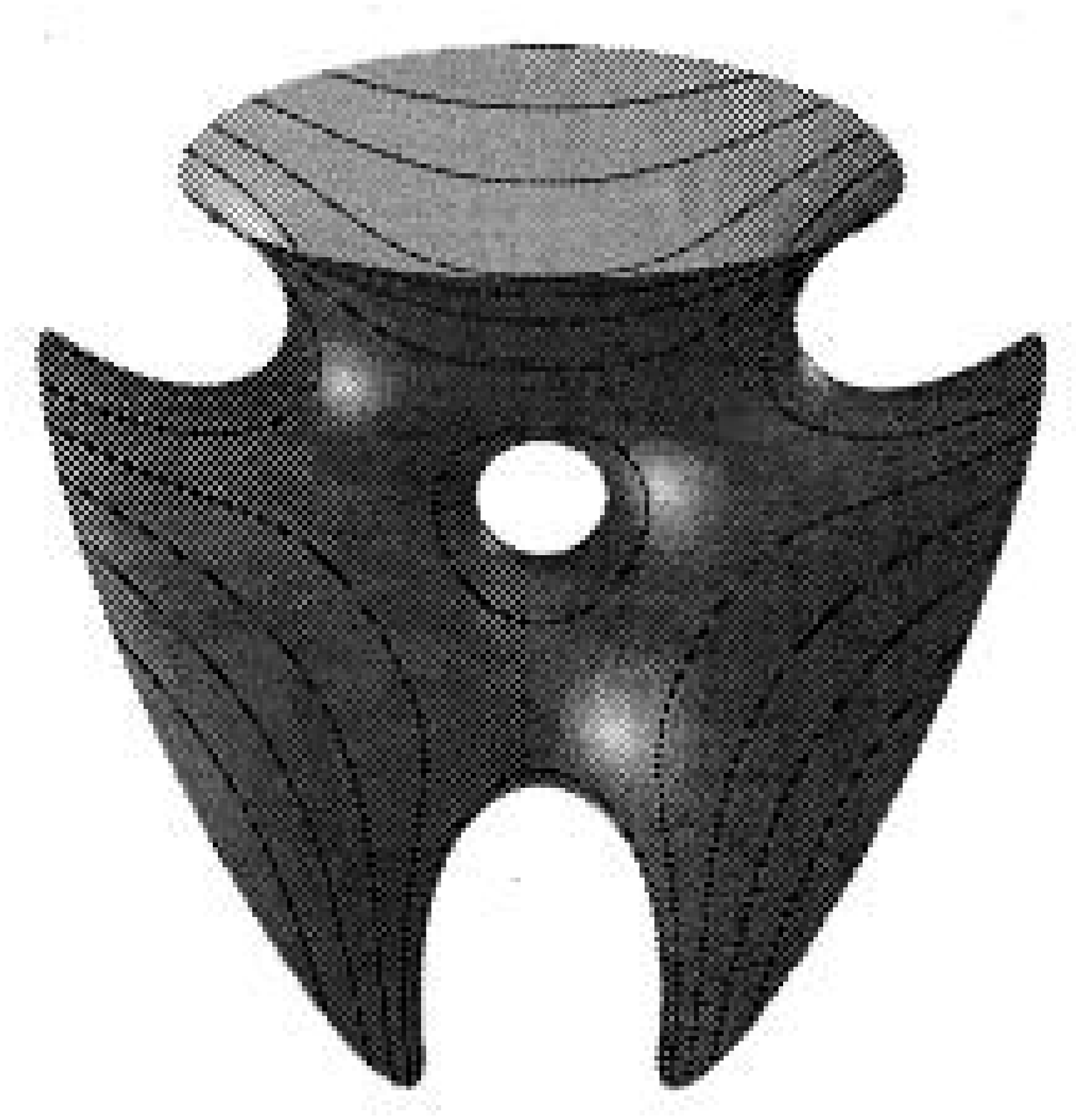}
        \hspace{0.01in}
        \epsfxsize=1.8in
        \epsffile{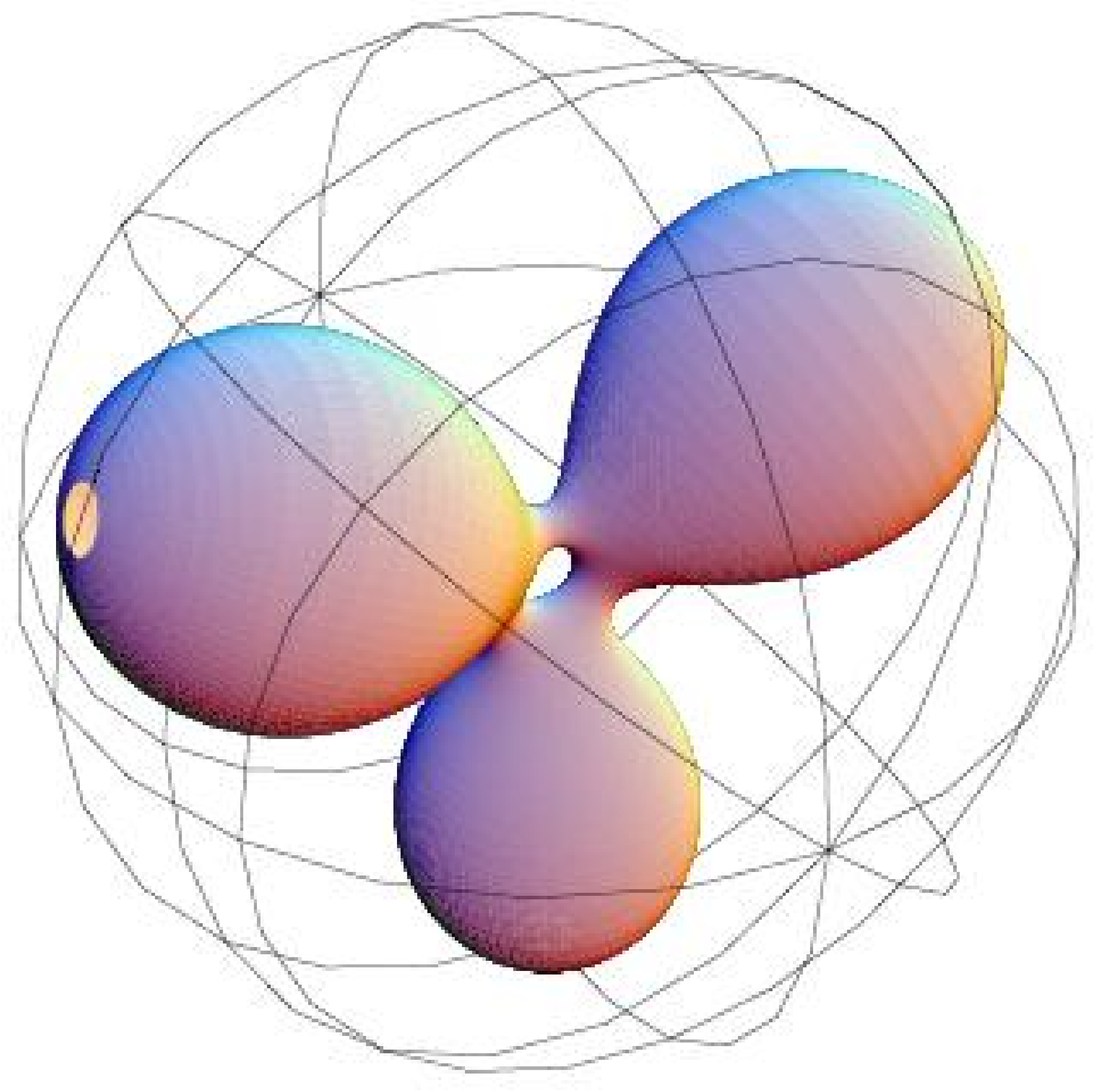}
        \hfill
\caption{A minimal genus 1 trinoid in $\bfR^3$, and a constant mean 
curvature 1 genus 1 trinoid cousin in $\bfH^3$.}
\end{figure}

\begin{figure}
        \hspace{0.5in}
        \epsfxsize=2.2in
        \epsffile{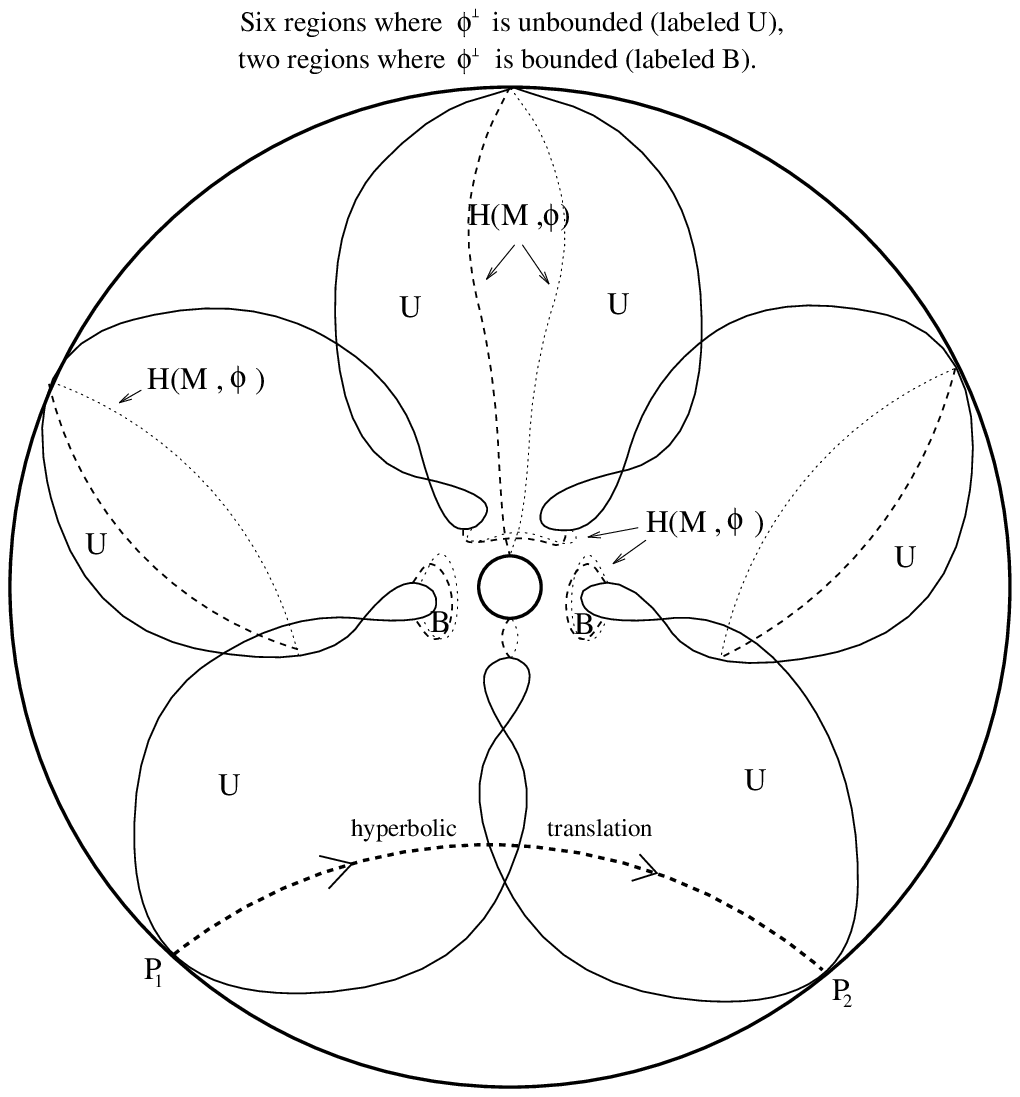}
        \hspace{0.5in}
        \epsfxsize=2.5in
        \epsffile{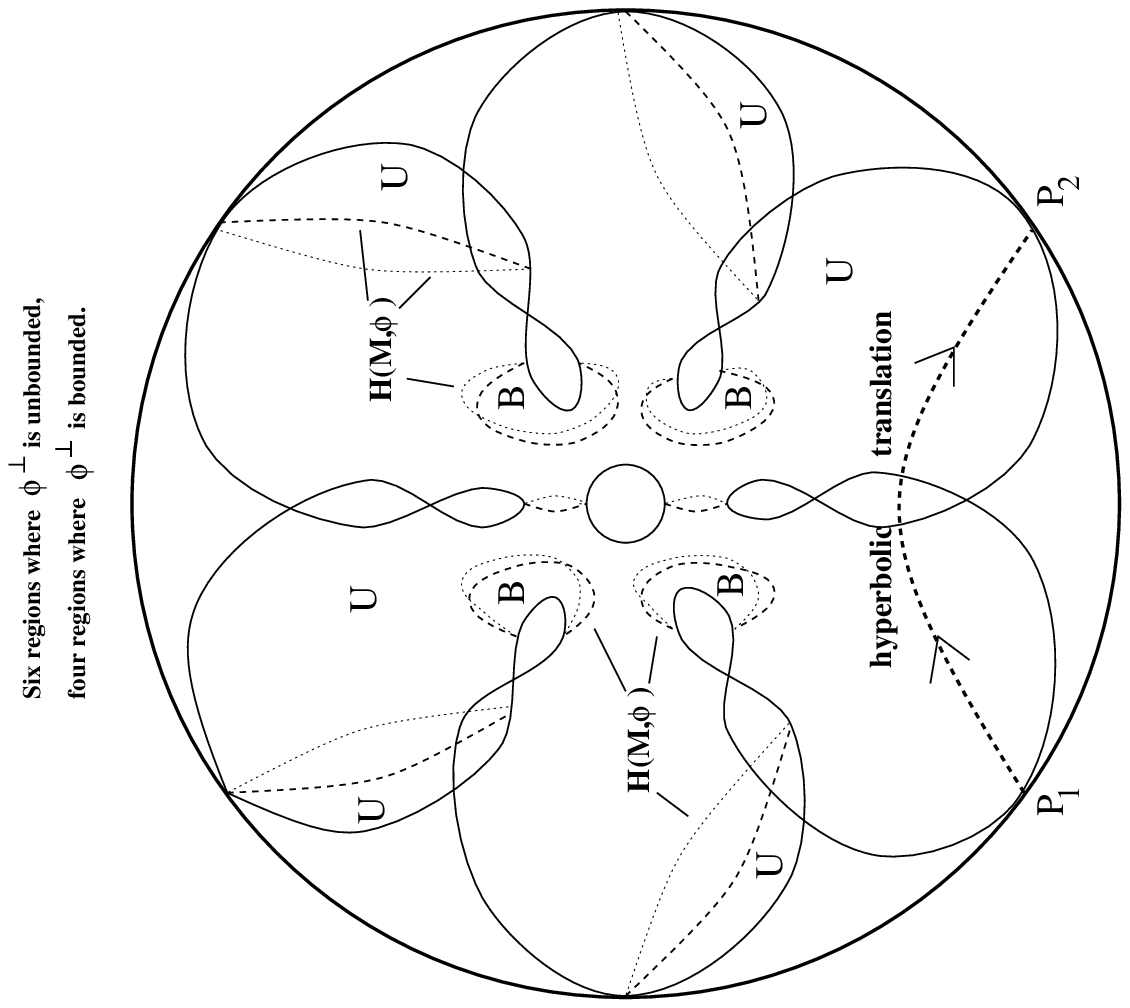}
        \hfill
\caption{The horizon $H(M,\phi)$ on genus 1 5-noid and 6-noid 
cousins, as described in the proof of Corollary 7.3.}
\end{figure}

\section{Deformations from minimal surfaces}

As stated in section 3, it was shown in \cite{RUY} that minimal 
surfaces in $\bfR^3$ can be deformed into corresponding 
constant mean curvature $c$ surface in $\bfH^3(-c^2)$.  And (as 
described in section 3) constant mean curvature $c$ surfaces in 
$\bfH^3(-c^2)$ are equivalent to constant mean curvature $1$ surfaces 
in $\bfH^3$.  Most of the known examples of complete 
constant mean curvature $1$ surfaces with finite topology 
in $\bfH^3$ have been shown to exist only via this deformation.  It is 
for this reason that the following theorem is of interest.  

For a minimal surface $M_0$, we consider the index Ind($M_0$) to be as 
defined in \cite{FC}.  The index of minimal surfaces is considered 
without a volume constraint, and this is natural because a volume 
constraint does not have a physical meaning for minimal surfaces.  

\begin{theorem}
If $M_0$ is a minimal surface in $\bfR^3$, and $M_c$ is
a corresponding constant mean curvature $c$ surface 
in $\bfH^3(-c^2)$, then if $c$ is
sufficiently close to zero, Ind($M_c$) $\geq$ Ind($M_0$)$-1$.  
\end{theorem}

It is not known yet if Ind($M_0$)$-1$ is the best possible lower 
bound for Ind($M_c$) in general, but 
the embedded catenoid cousins show that the best possible lower 
bound for Ind($M_c$) cannot be greater than Ind($M_0$) (since the 
index of a minimal catenoid is 1 and the index of an embedded catenoid 
cousin is also 1).  

Also, it is not possible to have 
an upper bound for the index of those constant mean curvature 
$c$ surfaces in
$\bfH^3(-c^2)$ which are deformations of minimal surfaces.  For 
example, the deformed Enneper cousin duals have infinite
total curvature for all $c \neq 0$, hence also infinite index 
(\cite{CS}), while the minimal Enneper cousin has index 1.

\begin{proof}
Let $\lambda_i(R,c)$ be the i'th eigenvalue of the operator 
$L$ with respect to the Dirichlet problem on $M_c \cap B_R(0)$, where 
$B_R(0) := \{(x_1,x_2,x_3) \in \bfR^3 \; | \; x_1^2 + x_2^2 + 
x_3^2 \leq R^2 \}$.  Note that $B_R(0)$ is contained in the Poincare 
model for $\bfH^3(-c^2)$ if $c$ is sufficiently close to zero, since 
the Poincare model has Euclidean radius $\frac{1}{c}$.  
Choose a finite 
$R>0$ large enough that Ind($M_{0}$) = Ind($M_0 \cap B_R(0)$) with 
Dirichlet boundary conditions.  
Thus \[ \lambda_1(R,0) \leq \dots \leq \lambda_k(R,0) < 0 \leq 
\lambda_{k+1}(R,0) \; . \]  Let $u_i = u_i(R,c)$ be the 
eigenfunction on $M_c \cap B_R(0)$ 
corresponding to the eigenvalue $\lambda_i(R,c)$.  

Note that for any finite value of $R$, $\lambda_i(R,c)$ is 
continuous in $c$.  We can see this
from the variational characterization of the eigenvalues (\cite{Be}, p60):
\[ \lambda_i(R,c) = \inf \{ \frac{\int_{M_c \cap B_R(0)} u Lu dA}
{\int_{M_c \cap B_R(0)} u^2 dA} 
\; | \; u \in C_0^\infty(M_c \cap B_R(0)), u \neq 0, u 
\perp_{L^{2}_{ds^2}} \{ u_{1},\ldots,u_{i-1} \} \; \} \; . \]
Since $M_c \cap B_R(0)$ is compact and $L$ is uniformly
continuous on a compact region, we know that this infimum will vary
continuously in $c$.
Therefore, for $c$ sufficiently close to 0, we still have 
\[ \lambda_1(R,c) \leq \dots \leq \lambda_k(R,c) < 0 \; . \]

Thus we have Ind$_u$($M_c$) $\geq$ Ind$_u$($M_c \cap B_R(0)$) $\geq$ 
Ind($M_0 \cap B_R(0)$) $=$ Ind($M_0$), and by Lemma 5.3 we have 
Ind($M_c$) $\geq$ Ind($M_0$)$-1$.  
\end{proof}

\begin{corollary}
For sufficiently small (constant mean curvature $c$) deformations $M$ 
in $\bfH^3(-c^2)$ from minimal surfaces in $\bfR^3$, we have 
the following lower bounds for index:
\begin{itemize}
\item If $M$ is a genus 0 $n$-noid cousin, then Ind($M$)$\geq 2n-4$.  
\item If $M$ is a genus 1 $n$-noid cousin, then Ind($M$)$\geq n-2$.  
\item If $M$ is a genus $k$ Costa cousin with $k \leq 37$, then 
Ind($M$)$\geq 2k+2$.  
\end{itemize}
\end{corollary}

At first, the second part of this lemma may seem like a stronger 
result than Corollary 7.3, but in Corollaries 8.1 and 
7.3, we do not know how large a deformation is 
possible.  It is possible that Corollary 7.3 
will allow larger deformations than Corollary 8.1 will allow.  
Hence we cannot say that Corollary 8.1 is a stronger 
result than Corollary 7.3.  

\begin{proof}
Nayatani \cite{N1} showed that the minimal genus 0 $n$-noid in 
$\bfR^3$ has index $2n-3$.  By Proposition 7.1, the minimal genus 1 
$n$-noid in $\bfR^3$ has index at 
least $n-1$.  Nayatani \cite{N2} also showed that the minimal genus $k$ 
Costa surface in $\bfR^3$ has index $2k+3$ for all $k \leq 37$.  
By Theorem 8.1, the corollary follows.  
\end{proof}

\section{Appendix: asymptotic behavior of ends}

Let $D \setminus \{0\}$ be the 
unit disk in the plane with the origin removed.  Let $\Phi: D \setminus 
\{0\} \rightarrow \bfH^{3}$ be a constant mean curvature 1, finite 
total curvature surface with a complete regular end 
at $0$.  We can take the secondary Gauss map to be $G = z^{\mu} 
\hat{G}$, and we can take the Hopf differential $Q$ so that 
$\frac{Q}{dG} = \omega dz = 
z^{\nu} \hat{\omega} dz$ where $\mu,\nu \in \bfR$ and $\hat{G}, 
\hat{\omega}$ are holomorphic and 
$\hat{G}(0) \neq 0, \hat{\omega}(0) \neq 0$ \cite{UY1} .  
As in the proof of Lemma 4.2, we may assume $\mu > 0$.  
The fact that the end is regular implies that ord$_0(Q) \geq -2$ 
\cite{By}, hence 
\[ Q = (\frac{q_{-2}}{z^2} + \frac{q_{-1}}{z} + \ldots) (dz)^2 \; . 
\] The leading coefficient $q_{-2}$ may or may not be zero.  
Completeness implies that $\nu \leq -1$ 
(this follows just by considering $ds^{2}$), and the fact that $Q$ is 
meromorphic implies that $\mu + \nu \in \bfZ$.  Finding a 
solution
\[ F^{-1} = \left( \begin{array}{cc} 
A & B \\ 
C & D 
\end{array} \right) \] in the Weierstrass representation 
(Lemma 3.2) for this surface, we have that 
$A$ and $C$ satifsy (see \cite{UY1}) \[ X^{\prime \prime} - 
\frac{\omega^{\prime}}{\omega} X^{\prime} - \omega G^{\prime} X 
 = 0 \; \; , \] (where $\prime$ denotes $\frac{d}{dz}$) 
and $B$ and $D$ satifsy \[ Y^{\prime \prime} - 
\frac{(G^{2}\omega)^{\prime}}{G^{2}\omega} Y^{\prime} - 
\omega G^{\prime} Y = 0 \; \; . \]  

The indicial equations of the above second order equations are 
\[ t^{2} - (\nu + 1) t - q_{-2} = 0  \mbox{  and  } 
t^{2} - (2 \mu + \nu + 1) t - q_{-2} = 0 \; \; . \]
The differences of solutions for the indicial equations are 
$m_{1} = \sqrt{(\nu + 1)^{2} + 4 q_{-2}}$ and $m_{2} = 
\sqrt{(2 \mu + \nu + 1)^{2} + 4 q_{-2}}$.  By results in 
\cite{UY1}, if the end is well-defined, then $m_{1},m_{2} \in 
\bfZ^+$, and the end 
is embedded if and only if $m:=$min($m_{1},m_{2}$) = 1.  

On page 626 of \cite{UY1} the end $\Phi$ is classified into 
three possible cases.  
\begin{itemize}
\item $\mu = 0$, $m_1 = m_2 = m$.  
\item $\mu \neq 0$, ord$_0 (Q) = -2$, $m_1 = m_2 = m$.  
\item $\mu \neq 0$, ord$_0 (Q) \geq -1$, $m_1 = -(\nu + 1)$, 
$m_2 = 2\mu + \nu + 1$, $m_2 - m_1 = 2($ord$_0 (Q) + 2) > 0$.  
\end{itemize}
In all three cases $m = m_1$.  

As we saw in the proof of Lemma 4.2, we may replace F by $BF$ for 
some $B \in$ SU(2) so that $\mu > 0$.  Thus the first case above 
always reduces to the second case (if ord$_0$($Q$)$=-2$) or the third 
case (if ord$_0$($Q$)$>-2$).  Hence we only need to consider cases 
2 and 3.  

In case 2, we say that we 
have a {\em catenoid cousin type end}.  In case 3 we say 
that we have a {\em horosphere type end}.  
Let ${\cal O}(z^\alpha)$ denote any 
complex valued function $f$ such that $\limsup_{z \to 0} 
\frac{f}{z^\alpha}$ is bounded.  

\begin{lemma}
At a catenoid cousin type end, 
$\mu \neq m$ and $F^{-1}$ can be locally represented as 
\[ F^{-1} = \frac{1}{\sqrt{\mu m}} \left( \begin{array}{cc} 
\frac{\mu+m}{2} z^{\frac{-\mu+m}{2}} (1+{\cal O}(z)) & 
\frac{\mu-m}{2} z^{\frac{\mu+m}{2}} (1+{\cal O}(z)) \\ 
\frac{\mu-m}{2} z^{\frac{-\mu-m}{2}} (1+{\cal O}(z)) & 
\frac{\mu+m}{2} z^{\frac{\mu-m}{2}} (1+{\cal O}(z)) 
\end{array} \right) \; \; . \]  
\end{lemma}

\begin{proof}
In case 2, $q_{-2} \neq 0$ and $\mu + \nu = -1$ and 
$m^2 = m_1^2 = (\nu + 1)^2 + 4 q_{-2}$.  Hence 
$(-\mu)^2 = (\nu+1)^2 \neq m^2$, so $\mu \neq m$.  

By Lemma 5.3 of \cite{UY1}, 
\[ \Delta \cdot F^{-1} = \frac{1}{\sqrt{\mu m}} 
\left( \begin{array}{cc} 
\frac{\mu+m}{2} z^{\frac{-\mu+m}{2}} a(z) & 
\frac{\mu-m}{2} z^{\frac{\mu+m}{2}} b(z) \\ 
\frac{\mu-m}{2} z^{\frac{-\mu-m}{2}} c(z) & 
\frac{\mu+m}{2} z^{\frac{\mu-m}{2}} d(z) 
\end{array} \right) \; , \] where $a(z),b(z),c(z),d(z)$ are 
holomorphic and nonzero at $z=0$.  Since $\Delta \in$ SL$(2,\bfC)$ 
only represents an 
isometry of $\bfH^{3}$ \cite{UY1}, we can assume $\Delta$ is the identity 
matrix.  By doing the transformation $z \rightarrow 
(G(0))^{\frac{-1}{\mu}} z$, we have $G = z^{\mu} (1+g_{1} z + 
\ldots)$.  By equation (3.1) we know 
$g = \frac{dA}{dC} = \frac{dB}{dD}$, and computing $g$ we find that 
$g \approx z^{m}$ and $a(0) d(0) = b(0) c(0)$.  
Since det($F^{-1}$) = 1, we 
have $a(0) d(0) = b(0) c(0) = 1$.  Since $G = \frac{-dB}{dA} = 
\frac{-dD}{dC}$, we have that $a(0) = b(0)$ and $c(0) = d(0)$.  
Doing an isometry of $\bfH^{3}$ so that 
\[ F^{-1} \rightarrow \left( \begin{array}{cc} 
\frac{1}{a(0)} & 0 \\ 0 & a(0) \end{array} \right) F^{-1} \; , \] we have 
proved the lemma.  
\end{proof}

\begin{lemma}
At a horosphere type end, $F^{-1}$ can be locally represented as 
\[ F^{-1} = \left( \begin{array}{cc} 
1+ {\cal O}(z) & {\cal O}(z^{2\mu+\nu+1}) \\ 
z^{\nu+1} (1  + {\cal O}(z)) & 1+ {\cal O}(z)
\end{array} \right) \; \; . \] 
\end{lemma}

\begin{proof}
Note that for case 3, we have 
$\mu \geq - \nu \geq 2$, $\mu,\nu \in \bfZ$.  
By Lemma 5.3 of \cite{UY1}, 
\[ \Delta \cdot F^{-1} = \left( \begin{array}{cc} 
a(z) & z^{2\mu+\nu+1} b(z) \\ 
z^{\nu+1} c(z) & d(z) 
\end{array} \right) \; , \] where $a(z),b(z),c(z),d(z)$ are 
holomorphic and nonzero at $z=0$.  Again, $\Delta$ only represents an 
isometry of $\bfH^3$.  We may assume $\Delta_{11} = \zeta \in 
\bfC \setminus \{0\}$, $\Delta_{22} = \zeta^{-1}$, $\Delta_{12} = 
\Delta_{21} = 0$, and then for any value of $\zeta$ the end is still 
asymptotic to the origin in the upper half space model.  With an 
appropriate choice of $\zeta$, we can conclude that $a(0) = d(0)$.  We 
can then rewrite $\Delta \cdot F^{-1}$ simply as $F^{-1}$.  
Since det($F^{-1}$) = 1, we have $a(0) = d(0) = 1$.  
Transforming $z \rightarrow 
(c(0))^{\frac{-1}{\nu+1}} z$, we have the lemma.  
\end{proof}

The point \[ F^{-1} \overline{F^{-1}}^t = \left( 
\begin{array}{cc} A & B \\ C & D 
\end{array} \right) \left( \begin{array}{cc} 
\bar{A} & \bar{C} \\ \bar{B} & \bar{D} 
\end{array} \right)  \] in the Hermitean model corresponds to the 
point 
\[ \frac{(\mbox{Re}(A \bar{C}+B \bar{D}),\mbox{Im}(A \bar{C}+B 
\bar{D}),
1)}{C \bar{C} + D \bar{D}}  \] in the upper half space model.  
So the catenoid cousin type end is 
\[ \frac{\mu+m}{\mu-m} \left(\mbox{Re}(z^{m}) 
(1+{\cal O}(|z|^{\mbox{min}(1,2\mu)})),
\mbox{Im}(z^{m}) (1+{\cal O}(|z|^{\mbox{min}(1,2\mu)})),
\frac{4 \mu m}{\mu^2-m^2} |z|^{\mu+m} 
(1+{\cal O}(|z|^{\mbox{min}(1,2\mu)})) \right) \] in the upper half space 
model. ${\cal O}(|z|^\alpha)$ denotes any real valued function such that 
$\limsup_{z \to 0} \frac{f}{|z|^\alpha}$ is finite.  The horosphere end 
is 
\[ \left(\mbox{Re}(z^{m}) (1+{\cal O}(|z|)),\mbox{Im}(z^{m}) 
(1+{\cal O}(|z|)),|z|^{2m} (1+{\cal O}(|z|)) \right) \] in the 
upper half space model, where $m = -\nu-1$.

\begin{lemma}
An end of the form 
\[\left(c_{1} \mbox{Re}(z^{m}) (1+{\cal O}(|z|^{\alpha})), 
c_{1} \mbox{Im}(z^{m}) (1+{\cal O}(|z|^{\alpha})), 
c_{2} |z|^{\mu+m} (1+{\cal O}(|z|^{\alpha})) \right)\] 
can be written in the form 
\[ \left(\mbox{Re}(z^{m}),\mbox{Im}(z^{m}), c_{3} |z|^{\mu+m} 
(1+{\cal O}(|z|^{\alpha})) \right) \; \; . \]
\end{lemma}

\begin{proof}
There exists $\tilde{z} \approx z$ so that 
$z^{m} (1+{\cal O}(|z|^{\alpha})) = \tilde{z}^{m}$, by the 
Weierstrass preparation theorem.  If follows that 
$(\frac{z}{\tilde{z}})^{m} - 1 \in {\cal O}(|z|^{\alpha})$, and 
therefore $\lim_{z \rightarrow 0} \frac{z}{\tilde{z}} = 
\lim_{z \rightarrow 0} (\frac{z}{\tilde{z}})^{\alpha} = 1$.  Also, we 
have $z^{m} - \tilde{z}^{m} \in {\cal O}(|z|^{\alpha}) z^{m}$, 
so $m z^{m-1} (z - \tilde{z}) \approx z^{m} - \tilde{z}^{m} \in 
{\cal O}(|z|^{\alpha + m})$ by the mean value theorem, 
so $z - \tilde{z} \in {\cal O}(|z|^{\alpha + 
1})$.  Now we have
\[ \left| \; c_{2} |z|^{\mu+m} (1+{\cal O}(|z|^{\alpha})) - 
c_{2} |\tilde{z}|^{\mu+m} \; \right| \leq c_{2} \left| \; |z|^{\mu+m} - 
|\tilde{z}|^{\mu+m} + |z|^{\mu+m} {\cal O}(|z|^{\alpha}) \; \right| \leq 
\] \[ 
c_{2} \left| \; |z|^{\mu+m} - |\tilde{z}|^{\mu+m} \; \right| + 
{\cal O}(|z|^{\alpha+\mu+m}) \leq 
\]\[ 
c_{2} \left| z^{\mu+m} - \tilde{z}^{\mu+m} \right| + 
{\cal O}(|z|^{\alpha+\mu+m}) \; \; \; \; \; \; \; \; \; \; \; \; \; 
\mbox{(triangle inequality)}
\] \[
\approx c_{2} (\mu+m) |z|^{\mu+m-1} |z-\tilde{z}| + 
{\cal O}(|z|^{\alpha+\mu+m}) = 
\]\[ 
|z|^{\mu+m-1} 
{\cal O}(|z|^{\alpha+1}) + {\cal O}(|z|^{\alpha+\mu+m}) 
= {\cal O}(|z|^{\alpha+\mu+m}) \; . \]
So we have that $c_{2} |z|^{\mu+m} (1+{\cal O}(|z|^{\alpha})) = 
c_{2} |\tilde{z}|^{\mu+m} + {\cal O}(|z|^{\alpha+\mu+m})$ and 
therefore $c_{2} |z|^{\mu+m} (1+{\cal O}(|z|^{\alpha})) = 
c_{2} |\tilde{z}|^{\mu+m}(1 + {\cal O}(|\tilde{z}|^{\alpha}))$.  
So we can rewrite our parametrization as 
\[ \left(c_{1} \mbox{Re}(\tilde{z}^{m}), c_{1} 
\mbox{Im}(\tilde{z}^{m}), c_{2} |\tilde{z}|^{\mu+m} 
(1 + {\cal O}(|\tilde{z}|^{\alpha})) \right) \; . \]
Then making the transformation $\tilde{z} = c_{1}^{-\frac{1}{m}} z$, 
we have finished the proof.   
\end{proof}

So a regular end of finite total curvature is of the form 
\[ \left( \mbox{Re}(z^{m}), \mbox{Im}(z^{m}), c |z|^{\mu+m} 
(1 + {\cal O}(|z|^{\mbox{min}(1,2\mu)})) \right) \; \; . \]  
Note that when we describe 
the end as a graph like this, the immersion is no longer conformal.  
The end is embedded if and only if
$m = 1$ \cite{UY1}.

\smallskip 

{\bf Addendum.} Lucas Barbosa and Pierre Berard have recently announced 
a result that would imply Ind($M$) = Ind$_u$($M$) for every case in 
this paper.  Using their result, one could strengthen Lemma 5.3 by 
excluding the Ind($M$) = Ind($\Sigma$)-1 case.  This would result in 
corresponding strengthenings of Theorem 6.1 (the $2 [ \mu ]$ case could 
be excluded), Theorem 7.1 ($\tilde v(M,\varphi) - 1$ and 
$\tilde v(M,\varphi) - 2$ could be replaced with $\tilde v(M,\varphi)$ and 
$\tilde v(M,\varphi) - 1$, respectively), Corollary 7.1 ($2k$ could be 
replaced with $2k+1$), Corollary 7.2 ($2$ could be replaced with $3$), 
Corollary 7.3 ($n-3$ and $n-4$ could be replaced with $n-2$ and $n-3$, 
respectively), Theorem 8.1 (Ind($M_0$)-1 could be replaced with 
Ind($M_0$)), and Corollary 8.1 ($2n-4$, $n-2$, and $2k+2$ could be 
replaced with $2n-3$, $n-1$, and $2k+3$, respectively).

\vspace{0.3in}

levi@mat.ufc.br

wayne@math.kobe-u.ac.jp

\end{document}